\documentclass[11pt,reqno]{amsart}
\usepackage{amscd
}%

\newcommand{\p}{\partial}
\newcommand{\Z}{\mathbb Z}
\newcommand{\Q}{\mathbb Q}
\renewcommand{\H}{\mathcal H}
\newcommand{\G}{\mathcal G}
\newcommand{\A}{\mathcal A}

\newcommand{\B}{\mathcal B}
\newcommand{\M}{\mathcal M}
\newcommand{\R}{\mathcal R}
\newcommand{\g}{\mathfrak g}
\newcommand{\yhom}{\tilde{Y}}
\newcommand{\ahom}{\tilde{A}}
\newcommand{\rhom}{\tilde{\rho}}
\renewcommand{\hbar}[1]{\tilde{\bar{#1}}}
\renewcommand{\O}{\operatorname O}
\newcommand{\im}{\operatorname{im}}
\newcommand{\ep}{\varepsilon}
\renewcommand{\phi}{\varphi}

\newcommand{\lk}{\operatorname{lk}}
\newcommand{\tr}{\operatorname{tr}}

\newcommand{\su}{\mathfrak{su}}
\newcommand{\so}{\mathfrak{so}}

\newcommand{\cs}{\mathbf{cs}}
\newcommand{\ad}{\operatorname{ad}}

\renewcommand{\u}{\mathfrak u}
\newcommand{\Hom}{\operatorname{Hom}}
\newcommand{\Hess}{\operatorname{Hess}}
\newcommand{\Herm}{\operatorname{Herm}}
\newcommand{\Stab}{\operatorname{Stab}}
\newcommand{\Sym}{\operatorname{Sym}}

\newcommand{\hol}{\operatorname{hol}}
\newcommand{\sign}{\operatorname{sign}}
\newcommand{\ZZ}{\Z_2\oplus \Z_2}
\newcommand{\PR}{\operatorname{PR}}
\renewcommand{\P}{\mathcal{PR}}
\newcommand{\Id}{\operatorname{Id}}%
\newtheorem{theorem}{Theorem}[section]
\newtheorem{lemma}[theorem]{Lemma}
\newtheorem{proposition}[theorem]{Proposition}
\newtheorem{corollary}[theorem]{Corollary}%
\theoremstyle{definition}
\newtheorem{remark}[theorem]{Remark}
\title{Rohlin's invariant and gauge theory, I.\\ Homology $3$-tori.}%
\thanks{The first author was partially supported by NSF Grants
9971802 and 0204386. The second author was partially supported by NSF
Grant 0196523.}
\author[Daniel Ruberman]{Daniel Ruberman}
\address{Department of Mathematics, MS 050\newline\indent Brandeis
University \newline\indent Waltham, MA 02454}
\email{\rm{ruberman@brandeis.edu}}
\author[Nikolai Saveliev]{Nikolai Saveliev}
\address{Department of Mathematics\newline\indent
University of Miami \newline\indent PO Box 249085
\newline\indent Coral Gables, FL 33124}
\email{\rm{saveliev@math.miami.edu}}


\begin{document}%

\begin{abstract}
This is the first in a series of papers exploring the relationship
between the Rohlin invariant and gauge theory. We discuss a
Casson-type invariant of a $3$-manifold $Y$ with the integral 
homology of the $3$-torus, given by counting projectively flat 
$U(2)$-connections. We show that its mod 2 evaluation is given by 
the triple cup product in cohomology, and so it coincides with a 
certain sum of Rohlin invariants of $Y$. Our counting argument 
makes use of a natural action of $H^1(Y;\Z_2)$ on the moduli space 
of projectively flat connections; along the way we construct 
perturbations that are equivariant with respect to this action.  
Combined with the Floer exact triangle, this gives a purely 
gauge-theoretic proof that Casson's homology sphere invariant 
reduces mod 2 to the Rohlin invariant.
\end{abstract}

\maketitle%

\section{Introduction}\label{S:intro}
Casson's introduction of his invariant for homology
3--spheres~\cite{akbulut-mccarthy,saveliev:casson} has had many
profound consequences in low-dimensional topology.  One of the
most important is the vanishing of the Rohlin invariant of a
homotopy sphere, which follows from Casson's identification of
his invariant, modulo $2$, with the Rohlin invariant of an
arbitrary homology sphere. The proof of this identification
proceeds via a surgery argument, in which a series of invariants
is defined for knots and links of several components. Ultimately,
these invariants are related to classical knot invariants, such
as the Alexander polynomial, and the theorem follows. This surgery
point of view, further developed in~\cite{walker:casson,lescop:casson},
finds its ultimate expression in the theory of finite-type invariants
of $3$-manifolds~\cite{le-murakami-ohtsuki,ohtsuki:finite}.

In this paper we give a proof of the equality of Rohlin's and
Casson's invariants (modulo $2$) in terms of the gauge theoretic
framework introduced by Taubes~\cite{taubes:casson}.  Many of the
ingredients for this proof are already in place, namely Taubes'
original work, and the surgery sequence of
Floer~\cite{floer:knots,braam-donaldson:knots} relating Casson-type
invariants of manifolds obtained by surgery on a knot.   Our main
contribution is to identify a Casson-type invariant of a homology
$3$-torus, $Y$, with a Rohlin-type invariant. This is accomplished by
relating the action of the group $H^1(Y;\Z_2)$ on the moduli space of
(projectively) flat connections to the cup product in cohomology.
The techniques we develop to deal with equivariant aspects of
non-smooth moduli spaces should be of independent interest in the
study of instanton Floer homology
(compare~\cite[\S5.6]{donaldson:floer}).

Let us briefly describe the invariants in question; more details will
be given in the next section.  By \emph{homology $3$-torus} we mean a
closed oriented $3$-manifold $Y$ having the integral homology of the
$3$-torus $T^3 = S^1\times S^1\times S^1$. For any non--trivial $w\in
H^2(Y;\Z_2)$, we consider projectively flat connections on a principal
$U(2)$ bundle $P\to Y$ whose associated $SO(3) = PU(2)$ bundle has
second Stiefel--Whitney class equal to $w$. We define a Casson-type
invariant $\lambda'''(Y,w)$ to be one-half of the signed count of such
connections, modulo an appropriate gauge group. This invariant is
one-half of the Euler characteristic of the Floer homology studied
in~\cite{donaldson:floer} and~\cite{braam-donaldson:knots} and is not,
\emph{a priori}, an integer.

A pair consisting of an oriented $3$-manifold $X$ and a spin structure
$\sigma$ has a Rohlin invariant $\rho(X,\sigma)\in \Q/2\Z$. By
definition,
$$
\rho(X,\sigma) = \frac{1}{8}\,\sign(V)
$$
for any spin $4$-manifold $V$ with (spin) boundary $(X,\sigma)$. By the
Rohlin invariant $\rho'''(Y)$ of a homology $3$-torus $Y$ we mean the
sum, over the eight spin structures on $Y$, of their Rohlin invariants.
It is easy to see that, as for a homology sphere, this invariant
actually takes values in $\Z/2\Z$.

\begin{theorem}\label{T:main}
For any choice of non-trivial $w\in H^2(Y;\Z_2)$, the Casson invariant
$\lambda'''(Y,w)$ is an integer. If $\{a_1, a_2, a_3\} $ is a basis
in $H^1(Y;\Z)$ then
$$\lambda'''(Y,w) = (a_1\cup a_2\cup a_3)\,[Y]\pmod 2.$$
\end{theorem}

Note that this implies that $\lambda'''(Y,w) \pmod 2$ is independent
of $w$, so long as $w$ is non-trivial. It is a theorem of V. Turaev~\cite{turaev:linking}, based on S. Kaplan
\cite[Lemma 6.3]{kaplan:even} that the above triple cup product also
evaluates the Rohlin invariant. Hence we obtain the following result.

\begin{corollary}\label{C:CRtori}
If $Y$ is a homology $3$-torus, then $\lambda'''(Y,w)\equiv\rho'''(Y)
\pmod 2$.
\end{corollary}

At the end of the paper, we will explain how this implies Casson's
original result about the Rohlin invariant of homology spheres.

We should point out that the congruence in Theorem~\ref{T:main}
actually holds over the integers (with appropriate sign conventions).
A proof of this stronger statement is implicit in Casson's original
work~\cite{akbulut-mccarthy,saveliev:casson} and a closely related
formula is given by Lescop~\cite{lescop:casson}.  To establish the
integral version of Theorem~\ref{T:main}, one would have to show the
equality of the count of flat connections on a homology $3$-torus
with either Lescop's invariant or Casson's invariant for 3-component
links (of trivial linking numbers). This equality is closely related
to Casson's formula for his knot invariant in terms of the Alexander
polynomial; a purely gauge-theoretic proof of the latter has been
given by Donaldson~\cite{donaldson:casson}. The techniques in the
present paper are rather different, and moreover have the advantage
of extending to the $4$-dimensional situation
\cite{ruberman-saveliev:donaldson} where the integral version does
not hold.

The idea of the proof of Theorem~\ref{T:main} is to take advantage of
a natural $H^1(Y;\Z_2)\allowbreak = (\Z_2)^3$ action on the moduli
space of projectively flat connections. We identify this moduli space
with the space of projective representations of $\pi_1(Y)$ in $SU(2)$,
and use this identification to show that the orbits with two elements
are always non--degenerate and that the number of such orbits equals
$(a_1\cup a_2\cup a_3)\,[Y]\pmod 2$. In the non--degenerate situation,
this completes the proof because there are no orbits with just one
element, and the orbits with four and eight elements do not contribute
to $\lambda'''(Y)\pmod 2$. The general case reduces to the
non--degenerate one after one finds a generic perturbation which
commutes with the $H^1(Y;\Z_2)$--action.  As mentioned above, this
equivariance is rather delicate.

The authors are thankful to Christopher Herald for sharing his
expertise.


\section{The invariant $\lambda'''$}

In this section we introduce the invariant $\lambda'''$ of a homology
3-torus $Y$ by counting projectively flat connections in a
$U(2)$--bundle over $Y$.  The `derivative' notation comes from
Casson's original approach, in which $\lambda'''$ appears as the
third difference quotient of his homology sphere invariant.


\subsection{The bundles}

Let $Y$ be a homology 3--torus, $P$ a principal $U(2)$ bundle over $Y$,
and $\bar P$ its associated $SO(3) = PU(2)$ bundle.  Topologically, the
bundles $P$ and $\bar P$ are determined by their characteristic classes
$c_1(P)$ and $w_2(\bar P)$, which are related by the formula $w_2(\bar
P) = c_1(P)\pmod 2$. Since $H^3(Y;\Z)$ is torsion free, every $SO(3)$
bundle over $Y$ arises as $\bar P$ for some $U(2)$ bundle $P$, and
$SO(3)$ bundles with non-trivial $w_2$ correspond to $U(2)$ bundles
whose $c_1$ is an odd element in $H^2(Y;\Z)$.

Every connection $A$ on $P$ induces connections on $\bar P$ and on the
line bundle $\det P$, via the splitting $\u(2) = \su(2)\,\oplus\,\u(1)$.
In a local trivialization, this corresponds to the decomposition
\begin{equation}\label{E:split}
A = \left(A - \frac 1 2\,\tr A\cdot\Id\right) + \frac 1 2\,\tr A\cdot\Id.
\end{equation}
The induced connection on $\bar P$ is the image of the first summand under
the isomorphism $\ad: \su(2) \to \so(3)$ given by $\ad(\xi)(\eta) = [\xi,
\eta]$, and the induced connection on $\det P$ is $\tr A$. Conversely, any
two connections on $\bar P$ and $\det P$ determine a unique connection on
$P$.

Fix a connection $C$ on $\det P$, and let $\A(P)$ be the space of connections
on $P$ compatible with $C$. The connection $C$ plays no real geometric role
-- different choices will give equivalent theories. The gauge group $\G(P)$
consisting of unitary automorphisms of $P$ of determinant one preserves $C$
and hence acts on $\A(P)$ with the quotient space $\B(P) = \A(P)/\G(P)$. Let
$\A(\bar P)$ be the affine space of connections on $\bar P$ and $\G(\bar P)$
the $SO(3)$ gauge group. Denote $\B(\bar P)= \A(\bar P)/\G(\bar P)$. The
projection
$\pi: \A(P) \to \A (\bar P)$ induced by the splitting~\eqref{E:split}
commutes with the above gauge group actions and hence defines a projection
\begin{equation}\label{E:pi}
\pi: \B(P) \to \B(\bar P).
\end{equation}

\medskip
The group $H^1(Y;\Z_2)$ acts on $\B(P)$ as follows. Let us view $\chi\in
H^1(Y;\Z_2)$ as a homomorphism from $\pi_1(Y)$ to $\Z_2 = \{\,\pm 1\,\}$.
As such, it defines a complex line bundle $L_{\chi}$.  For any $A \in
\A(P)$, let $A\otimes\chi$ be the connection on $P\otimes L_{\chi}$ induced
by $A$ and $\chi$. In a local trivialization, $A\otimes\chi$ is given by
$A + \omega$ where $\omega$ is an $\su(2)$-valued 1-form. In particular,
we easily see that both $A$ and $A\otimes\chi$ define the same connection
$C$ on $\det(P) = \det(P\otimes L_{\chi})$. Since the bundles $P\otimes
L_{\chi}$ and $P$ are isomorphic, the action is well defined on
gauge-equivalence classes by the formula $A\mapsto A\otimes\chi$.

\begin{proposition}
The map $\pi$ defined in \eqref{E:pi} is the quotient map of the
$H^1(Y;\Z_2)$ action described above.
\end{proposition}

\begin{proof}
The connections on $\bar P$ induced by $A$ and $A\otimes\chi$ are $SO(3)$
gauge equivalent because they have the same holonomy. Every connection on
$\bar P$ arises from a connection on $P$ which is unique up to the action
in question.
\end{proof}


\subsection{Projectively flat connections}\label{S:2.2}

Let $A$ be a connection on $P$ compatible with the connection $C$ on $\det
P$ and let $\bar A = \pi (A)$. The projection $\pi: \A(P) \to \A(\bar P)$
identifies the tangent spaces of $\A(P)$ and $\A(\bar P)$ at $A$ and $\bar A$,
respectively. The latter tangent space is known to be isomorphic to
$\Omega^1 (Y; \ad\bar P)$ where $\ad\bar P = \bar P \times_{\ad}\so(3)$. A
straightforward calculation shows that the curvatures of $A$ and $\bar A$
are related by

\[
F_A = \pi_*^{-1}(F_{\bar A}) + \frac 1 2\; F_C\cdot\Id.
\]

\medskip\noindent
We say that $A$ is a \emph{projectively flat connection} compatible with $C$
if $F_{\bar A} = 0$. The property of a connection being projectively flat is
preserved by the actions of both $\G(P)$ and $H^1(Y;\Z_2)$. The moduli space
of projectively flat connections will be denoted by $\M(P)$.

Fix a projectively flat connection $A_0$ compatible with $C$. Define the
\emph{Chern--Simons function} $\cs_{A_0}: \A(P)\to \mathbb R/\Z$ by the formula

\[
\cs_{A_0}(A) = \frac 1 {8\pi^2}\,\int_Y\;\tr\left(B\wedge dB + \frac 2 3\,
B\wedge B\wedge B\right),
\]

\medskip\noindent
where $B = A - A_0$ and $\tr: \su(2)\to\mathbb R$ is the trace function.
This function is invariant with respect to the $H^1(Y;\Z_2)$ action and
defines the Chern--Simons function $\cs_{A_0}: \B(\bar P) \to \mathbb R/4\Z$
on the quotient, so that we have the following commutative diagram

\[
\begin{CD}
\B(P)      @>\cs_{A_0} >> \mathbb R/\Z \\
@V\pi VV @V\times 4 VV \\
\B(\bar P) @>\cs_{A_0} >> \mathbb R/4\Z
\end{CD}
\]

\bigskip

The critical point set of $\cs_{A_0}: \B(P)\to\mathbb R/\Z$ can be
identified with the moduli space $\M(P)$ of projectively flat
connections on $P$, which is independent of the choice of $A_0$.
For this reason, we will generally omit the $A_0$ subscript in what
follows. The group $H^1(Y;\Z_2)$ acts on $\M(P)$. The quotient of
this action is the flat moduli space $\M(\bar P)$, which is the
critical point set of $\cs: \B(\bar P)\to \mathbb R/4\Z$. If
$w_2 (\bar P)\ne 0$ then all flat connections in $\M(\bar P)$ are
irreducible.


\subsection{Definition of $\lambda'''$}\label{S:def}

Let $A$ be a projectively flat connection on $P$. The point $[A]\in\M(P)$
is said to be \emph{non--degenerate} if $H^1(Y;\ad A) = 0$. Here,
$H^1 (Y;\ad A)$ stands for the cohomology with coefficients in the flat
bundle $\ad\bar P$ endowed with the flat connection $\bar A$. The moduli
space $\M(P)$ is called \emph{non--degenerate} if all of its points are
non-degenerate.

Let $w$ be a non-zero class in $H^2(Y;\Z_2)$, and $P$ a $U(2)$ bundle
over $Y$ with $c_1(P) = w\pmod 2$. If $\M(P)$ is non-degenerate then it
is finite and we define the \emph{Casson invariant} $\lambda'''(Y,w)$
as
\[
\lambda'''(Y,w) = \frac 1 2\,\sum_{A\in\M(P)}\; (-1)^{\mu(A)},
\]
where $\mu(A)$ is the mod 2 Floer index of $A$ defined as in
\cite[page 150]{donaldson:floer}. Note that the usual Floer index defined
modulo 8 is relative; for any pair of projectively flat connections $A_1$
and $A_2$, the modulo 2 reduction of this relative index equals $\mu(A_1)
- \mu(A_2)\pmod 2$.

If $\M(P)$ happens to be degenerate then it will need to be perturbed as
described in Section \ref{S:perturb}, and then $\lambda'''(Y)$ will be
defined essentially as above. That $\lambda'''(Y)$ is well defined follows
from \cite[pages 148--149]{donaldson:floer}.

\medskip

\begin{proposition}\label{P:sf}
The action of $H^1(Y;\Z_2)$ preserves the mod 2 Floer index.
\end{proposition}

\begin{proof}
This follows from \cite[pages 239--240]{braam-donaldson:knots}.
\end{proof}

\begin{remark}
According to Proposition~\ref{P:sf}, the points in the
$H^1(Y;\Z_2)$--orbit of a projectively flat connection $A$ are counted
in $\lambda'''(Y,w)$ with the same sign. Hence we could as well define
$\lambda'''(Y,w)$ by counting points in $\M(\bar P)$, where $w_2(\bar P)
= w$, with weights given by the order of the orbits of their respective
lifts to $\M(P)$.
\end{remark}


\section{Projective representations}\label{S:proj}

The holonomy map gives a homeomorphism between the moduli space $\M(\bar P)$
of flat connections on $\bar P$ and the $SO(3)$--character variety of
$\pi_1(Y)$.   Similarly, there is an algebraic interpretation (again using
holonomy) of projectively flat connections in terms of \emph{projective}
representations. This section describes this concept in some detail; good
general references for these ideas are the classic paper of
Atiyah-Bott~\cite{atiyah-bott:surface} and the book of Brown~\cite{brown}.


\subsection{Algebraic background}

Let $G$ be a finitely presented group and view $\Z_2 = \{\pm 1\}$ as the
center of $SU(2)$. A map $\rho: G\to SU(2)$ is called a \emph{projective
representation} if $\rho(gh)\rho(h)^{-1}\rho(g)^{-1}\in\Z_2$ for all $g, h
\in G$. Given a projective representation $\rho$, the function $c: G\times
G\to\Z_2$ defined as $c(g,h) = \rho(gh)\rho(h)^{-1}\rho(g)^{-1}$ is a
2--cocycle, that is, $c(gh,k)c(g,h) = c(g,hk)c(h,k)$. We will refer to $c$
as the \emph{cocycle associated with $\rho$}.

Let us fix a cocycle $c: G\times G\to\Z_2$ and denote by $\PR_c (G;SU(2))$
the set of all projective representations $\rho: G\to SU(2)$ whose
associated 2--cocycle is $c$.

\begin{lemma}\label{lemma1}
If $c$ and $c': G\times G\to \Z_2$ are cocycles such that $[c] = [c']\in
H^2(G;\Z_2)$ then there is a bijection between $\PR_c (G;SU(2))$ and
$\PR_{c'} (G;SU(2))$.
\end{lemma}

\begin{proof}
The fact that $[c] = [c']$ means that there exists a function $\mu: G\to
\Z_2$ such that $\mu(gh)c(g,h) = \mu(g)\mu(h)c'(g,h)$ for all $g, h\in G$.
Define a map $\phi: \PR_c (G;SU(2))\to \PR_{c'} (G;SU(2))$ by the formula
$\phi(\rho)(g) = \mu(g)\rho(g)$. One can easily check that $\phi(\rho)\in
\PR_{c'} (G;SU(2))$ and that $\phi$ is a bijection\,: its inverse $\psi:
\PR_{c'} (G;SU(2)) \to \PR_c (G;SU(2))$ is given by the formula
$\psi(\rho')(g) = \mu(g)\rho'(g)$.
\end{proof}

Let $c: G\times G\to\Z_2$ be a 2--cocycle and $\rho_1, \rho_2 \in
\PR_c(G;SU(2))$. We say that $\rho_1 \simeq \rho_2$ if there exists a
function $\mu: G\to\Z_2$ and an element $\sigma\in SU(2)$ such that
$\rho_2(g) = \mu(g)\sigma\rho_1(g)\sigma^{-1}$ for all $g\in G$.

\begin{lemma}\label{lemma3}
The map $\mu: G\to\Z_2$ is a homomorphism.
\end{lemma}

\begin{proof}
For any elements $g, h\in G$ we have $\rho_2(gh) = c(g,h)\rho_2(g)\rho_2(h)$.
This implies that $\mu(gh)\sigma\rho_1(gh)\sigma^{-1} = c(g,h)\mu(g)\sigma
\rho_1(g)\sigma^{-1}\mu(h)\sigma\rho_1(h)\sigma^{-1}$, and, after
simplification, $\mu(gh)\rho_1(gh) = \mu(g)\mu(h)c(g,h)\rho_1(g)\rho_1(h)$.
Since $\rho_1(gh) = c(g,h)\rho_1(g)\rho_1(h)$, we conclude that $\mu(gh) =
\mu(g)\mu(h)$.
\end{proof}

Let $\P_c(G;SU(2))$ be the set of conjugacy classes of projective
representations of $G$ viewed as $SU(2)$ valued functions. One can easily
see that $\simeq$ descends to equivalence relation on $\P_c(G;SU(2))$, and
hence there is a natural projection map
\begin{equation}\label{one}
\P_c (G;SU(2)) \quad\longrightarrow\;\P_c (G;SU(2))/\simeq \quad=\quad \PR_c
(G;SU(2))/\simeq.
\end{equation}

For any choice of cocycle $c: G\times G\to\Z_2$, the set $\PR_c (G;SU(2))$
is acted upon by the group $H^1(G;\Z_2) = \Hom(G;\Z_2)$. Every $\chi \in
\Hom(G;\Z_2)$ acts by the formula $\rho\mapsto\rho^{\chi}$ where
$\rho^{\chi}(g) = \chi(g)\rho(g)$, $g\in G$ (one can easily see that the
cocycle associated with $\rho^{\chi}$ is again $c$). This action preserves
conjugacy, and hence defines an action on $\P_c (G;SU(2))$.
\begin{proposition}
The quotient of $\P_c (G;SU(2))$ by the $H^1(G;\Z_2)$ action equals
$\P_c (G;SU(2))/\simeq$.
\end{proposition}

\begin{proof}
If $\rho_1\simeq\rho_2$ then there exists a map $\mu: G\to\Z_2$ and
an element $\sigma\in SU(2)$ such that $\rho_2(g) = \mu(g)\sigma
\rho_1(g)\sigma^{-1}$. Since $\mu$ is necessarily a homomorphism by
Lemma \ref{lemma3}, the above equality means that $\rho_2$ is
conjugate to $\rho_1^{\,\mu}$. The same formula shows that if $\rho_2$
is conjugate to $\rho_1^{\chi}$ for $\chi\in H^1(Y;\Z_2)$, then
$\rho_1\simeq\rho_2$.
\end{proof}

We now want to relate the projective $SU(2)$ representations studied
above to the ordinary $SO(3)$--representations of $G$. Let $\alpha:
G \to SO(3)$ be a representation. Denoting by $SO(3)^{\delta}$ the
group $SO(3)$ with the discrete topology, we have maps
$$
G\to SO(3)^\delta \to SO(3)
$$
induced by $\alpha$ and the identity map respectively. This yields maps
of classifying spaces
$$
BG\to BSO(3)^{\delta} \to BSO(3).
$$
Since $G$ is a discrete group, we can identify $H^2(BG;\Z_2)$ with the group
cohomology $H^2(G;\Z_2)$.  Thus we obtain a homomorphism
$$H^2(BSO(3);\Z_2)\to H^2(BSO(3)^\delta;\Z_2)\to
H^2(G;\Z_2).
$$
Let $w_2(\alpha)$ be the image in $H^2(G;\Z_2)$ of the universal
Stiefel--Whitney class $w_2\in H^2(BSO(3);\Z_2)$.
%
%
%

\begin{proposition}\label{P:w2}
Let $\ad\rho: G\to SO(3)$ be the composition of $\rho\in\PR_c (G;SU(2))$
and $\ad: SU(2)\to SO(3)$. Then $\ad\rho$ is a representation, and
$w_2 (\ad\rho) = [c]\in H^2(G;\Z_2)$.
\end{proposition}

\begin{proof}
This follows from the description of $w_2(\ad\rho)$ as the obstruction to
lifting $\ad\rho$ to an $SU(2)$ representation.
\end{proof}

\begin{corollary}
Let $\rho: G\to SU(2)$ be a projective representation with associated
2--cocycle $c$. Suppose that there is a non--central element $u\in SU(2)$
such that $u\rho(g) = \rho(g)u$ for all $g\in G$. Then $[c] = 0$ in $H^2
(G;\Z_2)$.
\end{corollary}

\begin{proof}
The image of $\rho$ is contained in a circle in $SU(2)$ hence $\ad\rho$ is
conjugate to an $SO(2)$ representation and hence admits an $SU(2)$ lift.
This implies that $[c] = w_2 (\ad\rho) = 0$.
\end{proof}

A projective representation $\rho: G\to SU(2)$ is called \emph{irreducible}
if the centralizer of its image equals the center of $SU(2)$. According to
the above corollary, any projective representation whose 2--cocycle is not
cohomologous to zero is irreducible.

Let $w\in H^2(G;\Z_2)$ and denote by $\R_w (G;SO(3))$ the set of the
conjugacy classes of $SO(3)$ representations of $G$ whose second
Stiefel--Whitney class equals $w$. This is a compact real algebraic
variety. The correspondence $\rho\mapsto\ad\rho$ defines a map
\begin{equation}\label{two}
\PR_c (G;SU(2))/\simeq\quad\longrightarrow\; \R_{[c]} (G;SO(3)).
\end{equation}

\begin{proposition}
The map (\ref{two}) is a bijection.
\end{proposition}

\begin{proof}
Suppose that $\rho_1,\rho_2\in\PR_c (G;SU(2))$ are such that $\ad\rho_1$
and $\ad\rho_2$ are conjugate as $SO(3)$ representations. Then there exists
a function $\mu: G\to \Z_2$ and and element $\sigma\in SU(2)$ such that
$\rho_2(g) = \mu(g)\sigma\rho_1(g)\sigma^{-1}$ for all $g\in G$. This
means that $\rho_1\simeq\rho_2$ and the map (\ref{two}) is injective.

Given a representation $\ad\rho: G\to SO(3)$, we can always lift it to
a projective representation $\rho'\in\PR_{c'} (G;SU(2))$ for some $c'$
such that $[c'] = [c]$. But then we can also find a lift $\rho\in
\PR_c (G;SU(2))$ because $\PR_{c'} (G;SU(2)) = \PR_c (G;SU(2))$ by Lemma
\ref{lemma1}.
\end{proof}


\subsection{The holonomy correspondence}\label{S:hol}

In this section we establish a correspondence between projectively flat
connections over a manifold $Y$ (which is not necessarily a homology
$3$-torus) and projective representations of its fundamental group. The
correspondence is, in rough terms, given by taking the holonomy of a
projectively flat connection. In principle, this is well-known, but we
could not find a reference. Moreover, some subtle points arise in
establishing the continuity of the correspondence.

In what follows, we will use the principle that connections pull back
under smooth maps. More precisely, let $j: M \to W$ be a smooth map
and suppose that $Q\to W$ is a principal $G$--bundle with connection,
determined by a $1$-form $\omega \in \Omega^1 (Q;\g)$, where $\g$ is
the Lie algebra of $G$. There is a bundle map $\bar j: j^*Q \to Q$
which commutes with the $G$ actions and which is an isomorphism on the
fibers. Then $\bar j^*\omega$ gives a connection on the bundle $j^*Q$,
whose holonomy has the following property: If $\gamma:I \to M$ is a
loop, then
$$
\hol_{j^*\omega}(\gamma) = \hol_{\omega}(j_*\gamma).
$$

The natural map $Y \to B(\pi_1 Y)$ induces a monomorphism $\iota:
H^2(\pi_1 Y; \Z_2) \to H^2 (Y;\Z_2)$, see \cite{brown}. We first
deal with the discrepancy arising from the fact that $\iota$ need
not be surjective.

\begin{lemma}\label{L:empty}
Let $P$ be a $U(2)$--bundle over a manifold $Y$ such that
$w_2(\bar P)$ is not in the image of $\iota: H^2 (\pi_1 Y;\Z_2)\to
H^2 (Y;\Z_2)$. Then the moduli space $\M(P)$ is empty.
\end{lemma}

\begin{proof}
The Hopf exact sequence $\pi_2 Y \to H_2 (Y;\Z) \to H_2 (\pi_1 Y;\Z)
\to 0$, see \cite{brown}, implies that, if $w_2 (\bar P)$ does not
belong to the image of $\iota$, it evaluates non-trivially on a
2-sphere in $Y$. Such a bundle $\bar P$ cannot support any flat
connections, for a flat connection on $\bar P$ would pull back to a
flat connection on the $2$--sphere, whose holonomy would trivialize
the bundle.
\end{proof}

  From now on, we will concentrate on bundles $P$ such that $w_2(\bar P)$
is in the image of $H^2 (\pi_1 Y;\Z_2)$, and will identify $H^2(\pi_1 Y;
\Z_2)$ with its (monomorphic) image in $H^2(Y;Z/2)$.

It is a well known fact that the holonomy defines a bijection $\bar\phi:
\M(\bar P) \to \R_w (Y; SO(3))$ where $w_2 (\bar P) = w \pmod 2$. Given a
$U(2)$--bundle $P$ with $c_1 (P) = w \pmod 2$, our immediate goal will be
to define an $H^1(Y;\Z_2)$-equivariant map $\phi: \M(P) \to \P_c (Y;
SU(2))$, where $[c] = w$, which makes the following diagram commute

\[
\begin{CD}
\M (P) @>\phi >> \P_c (Y; SU(2)) \\
@VV\pi V    @VV\pi V \\
\M(\bar P) @>\bar\phi >> \R_w (Y; SO(3))
\end{CD}
\]

\bigskip\noindent
Here, $\pi: \P_c(Y;SU(2)) \to \R_w(Y;SO(3))$ is the map (\ref{one})
followed by the bijection (\ref{two}), see Section \ref{S:proj}.
It is straightforward to define a map $A \to \phi(A)$ by lifting
$\bar\phi(\bar A)$ to a projective representation.  However, such an
assignment might not be continuous, because the choice of lifting is
not canonical.

Let $A$ be a projectively flat connection on $P$ whose central part is a
fixed connection $C$ on the linear bundle $\det P$. For any based loop
$\gamma$ in $Y$, we let
\begin{equation}\label{E:hol-map}
\phi(A)(\gamma) = \hol_A(\gamma)\cdot \hol_C(\gamma)^{-1/2}\in SU(2)
\end{equation}
where the square root of $\hol_C(\gamma)\in U(1)$ in the second factor
is defined as follows.
Let $\Omega (Y)$ be the monoid of based loops in $Y$, and fix a
representative $\gamma$ in each connected component $\Omega_{[\gamma]}$
so that
\[
\Omega (Y) = \bigsqcup_{[\gamma] \in \pi_1 (Y)}\; \Omega_{[\gamma]}.
\]
Choose a square root $\hol_C(\gamma)^{1/2}\in U(1)$ for each of the
$\gamma$ and define a (based) map
\[
h: (\Omega_{[\gamma]},\gamma) \to (U(1),1)
\]
by the formula $h(\alpha) = \hol_C(\alpha) \cdot \hol_C(\gamma)^{-1}$. If
$\pi: U(1)\to U(1)$ is the squaring map then we want to lift $h$ to
$\tilde h$ such that $h = \pi \circ \tilde h$ (given such a lift, we get
a square root of $\hol_C(\alpha)$ by the formula $\tilde h(\alpha) \cdot
\hol_C (\gamma)^{1/2}$). The obstruction to the above lifting problem is
given by
\begin{alignat}{2}
{\mathcal O} &\in H^1 (\Omega_{[\gamma]},\gamma; \Z_2)     & & \notag \\
      &= \Hom(H_1(\Omega_{[\gamma]},\gamma; \Z); \Z_2) & & \notag \\
      &= \Hom(\pi_1(\Omega_{[\gamma]},\gamma); \Z_2)  & & \notag \\
      &= \Hom(\pi_1(\Omega_*,*); \Z_2),               & &
         \quad\text{where $*$ is the trivial loop,} \notag \\
      &= \Hom(\pi_2(Y,*);\Z_2). & & \notag
\end{alignat}
In particular, we immediately see that this obstruction vanishes as long
as $\pi_2 (Y) = 0$.

\begin{lemma}
For any $U(2)$ bundle $P$ such that $c_1 (P) = w \ne 0 \pmod 2$, the
obstruction ${\mathcal O}$ is zero.
\end{lemma}

\begin{proof}
Without loss of generality we may assume that $\gamma = *$, the trivial
loop. The obstruction $\mathcal O$ can be described as follows. Given a
homotopy $\alpha_t$, $0\le t\le 1$, such that $\alpha_0 = \alpha_1 = *$
define $\tilde h(\alpha_0)=1$, and $\tilde h(\alpha_t)$ by path
lifting. Then ${\mathcal O}(\sigma) = \tilde h(\alpha_1) \in \Z_2$.
The image of ${\mathcal O}$ in $\Hom(\pi_2(Y,*);\Z_2)$ is gotten by
viewing a $2$-sphere $\sigma$ in $Y$ as such a path.

Now, given a class $\sigma \in \pi_2(Y)$, the obstruction to
extracting a root of the bundle $\det P$ is given by evaluation of
$\sigma^*(w_2(\det P))$ on $S^2$.  Since $w_2(\det P) = w_2(\bar P)$
and $\bar P$ admits a flat connection, the latter evaluation has to
be zero. Thus the bundle $\det P$ admits a square root over every
2-sphere, and the holonomy of this square root along the
loops $\alpha_t$ gives a lift $\tilde h$ which necessarily satisfies
$\tilde h(\alpha_0) = \tilde h(\alpha_1)$. In particular,
$\mathcal O$ vanishes.
\end{proof}

With the above definition of $\hol_C(\gamma)^{1/2}$ in place, the map $\phi$
is defined by the formula (\ref{E:hol-map}).  For any projectively flat
connection $A$ in $P$, composition of $\phi$ with the natural projection
$SU(2) \to SO(3)$ gives a representation $\bar\phi (A): \pi_1 (Y)\to
SO(3)$. In particular, if $[\gamma_1] = [\gamma_2]$ then $\phi (A)(\gamma_1)
= \pm \phi (A)(\gamma_2)$. Together with continuity of $\phi (A)$ this
implies that $\phi(A)(\gamma_1) = \phi(A)(\gamma_2)$ and hence $\phi(A):
\pi_1 (Y) \to SU(2)$ is a well defined projective representation.

\begin{proposition}\label{P:hol}
Let $P$ be a $U(2)$ bundle over $Y$ such that $c_1(P) = w\ne 0\pmod 2$. Then
$\phi$ is well defined as a map $\phi: \M(P) \to \P_c (Y; SU(2))$ with $[c] =
w$, and it is an $H^1 (Y; \Z_2)$--equivariant bijection.
\end{proposition}

\begin{proof}
The connection $C$ on $\det P$ does not change when $A$ is replaced by a
gauge equivalent connection. Therefore, the second factor in (\ref{E:hol-map})
remains unchanged and the first one changes by conjugation. Since $\hol_C
(\gamma)^{-1/2}$ is central, the entire $\phi(A)$ changes by conjugation.
Therefore, the map $\phi: \M(P) \to \P_c (Y; SU(2))$ is well defined. The
cocycle $c$ is determined by the choice of $\hol_C(\gamma)^{1/2}$ for the
representative loops $\gamma$, different choices leading to cohomologous
cocycles. That $[c] = w$ can be read off the definition of $\phi$.

Let $\chi \in H^1 (Y; \Z_2)$ and replace $A$ by $A \otimes \chi$ then $A$ and
$A\otimes \chi$ induce the same connection $C$ on $\det P$ so that the second
factor in (\ref{E:hol-map}) stays the same. The first factor becomes
$\hol_{A\otimes\chi}(\gamma) = \hol_A(\gamma)\cdot\chi(\gamma)$ with
$\chi(\gamma) = \pm 1$. Hence $\phi$ is $H^1(Y;\Z_2)$--equivariant. Since its
quotient map $\bar\phi$ is a bijection, so is $\phi$.
\end{proof}

An argument similar to that for representation varieties shows that Zariski
tangent space to $\P_c(Y;SU(2))$ at a projective representation $\rho:
\pi_1(Y)\to SU(2)$ equals $H^1 (Y;\ad\rho)$ where $\ad\rho: \pi_1(Y)\to SU(2)
\to SO(3)$ is a representation. It is identified as usual with the tangent
space to $\M(P)$ at the corresponding projectively flat connection.


\subsection{Application to homology $3$-tori}

A homology torus $Y$ is called \emph{odd} if there exist vectors $a_1$,
$a_2$, and $a_3$ in $H^1 (Y;\Z_2)$ such that $(a_1\cup a_2\cup a_3)\,[Y]
=1 \pmod 2$. Note that such $a_1$, $a_2$, and $a_3$ form a basis of
$H^1 (Y;\Z_2)$ because they are distinguished by cup-products with
$a_1\cup a_2$, $a_2\cup a_3$, and $a_1\cup a_3$ and hence are linearly
independent. Also note that if $(a_1\cup a_2\cup a_3)\,[Y] = 1 \pmod 2$
for some basis $a_1$, $a_2$, $a_3\in H^1 (Y;\Z_2)$ then the same is true
for any other basis. A homology torus $Y$ is called \emph{even} if
$(a_1\cup a_2\cup a_3)\,[Y] = 0 \pmod 2$ for any three vectors $a_1$,
$a_2$, $a_3\in H^1 (Y;\Z_2)$.

Let $\Lambda^2 H^1(Y;\Z_2)$ be the second exterior power of $H^1(Y;\Z_2)$
and consider the cup--product map
\begin{equation}\label{E:cup}
\cup\,: \Lambda^2 H^1(Y;\Z_2)\to H^2(Y;\Z_2).
\end{equation}

\begin{lemma}\label{L:1}
The map (\ref{E:cup}) is an isomorphism if $Y$ is odd, and is zero if $Y$
is even.
\end{lemma}

\begin{proof}
Let $Y$ be an odd homology torus and choose a basis $a_1$, $a_2$, $a_3\in
H^1 (Y;\Z_2)$. The vectors $a_1\cup a_2$, $a_2\cup a_3$, and $a_1\cup a_3
\in H^2 (Y;\Z_2)$ are linearly independent because they are distinguished
by the homomorphisms $H^2 (Y;\Z_2)\to H^3 (Y;\Z_2)$ given by cup-products
with $a_1$, $a_2$, and $a_3$. Suppose now that $Y$ is even and that there
are vectors $a$, $b\in H^1 (Y;\Z_2)$ such that $a\,\cup\, b\ne 0\pmod 2$.
By Poincar\'e duality, there exists $c\in H^1 (Y;\Z_2)$ such that $a\cup
b\cup c = 1\pmod 2$, a contradiction.
\end{proof}

\begin{corollary}\label{L:oddiota}
If $Y$ is odd then the map $\iota: H^2 (\pi_1 Y;\Z_2)\to H^2 (Y;\Z_2)$ is
an isomorphism.
\end{corollary}

\begin{proof}
This follows from the commutative diagram

\[
\begin{CD}
\Lambda^2 H^1 (\pi_1 Y;\Z_2) @>\cong >> \Lambda^2 H^1 (Y;\Z_2) \\
@V \cup VV      @V \cup VV \\
H^2 (\pi_1 Y;\Z_2) @> \iota >> H^2 (Y;\Z_2)
\end{CD}
\]

\medskip\noindent
whose upper arrow is an isomorphism because $H^1(\pi_1 Y;\Z_2) = H^1(Y;\Z_2)$,
and whose right arrow is an isomorphism by Lemma \ref{L:1}. Since $\iota$ is
injective, the remaining two arrows in the diagram are also isomorphisms.
\end{proof}

Note that if $Y$ is an even homology torus, the conclusion of Corollary
\ref{L:oddiota} need no longer hold\,: take for example $Y = (S^1\times S^2)
\# (S^1\times S^2)\# (S^1\times S^2)$.

\medskip

\begin{corollary}\label{C:iota}
Theorem \ref{T:main} holds for all $\lambda'''(Y,w)$ such that $w$ is not in
the image of $\iota: H^2 (\pi_1 Y;\Z_2)\to H^2 (Y;\Z_2)$.
\end{corollary}

\begin{proof}
Let $P$ be a bundle with $w_2(\bar P) = w$ not in the image of $\iota$.
Then, according to Lemma \ref{L:empty}, the moduli space $\M(P)$ is empty
so that $\lambda'''(Y,w) = 0$. On the other hand, this situation is only
possible if $Y$ is an even homology torus, see Corollary \ref{L:oddiota}.
\end{proof}

We will assume from now on that $w = w_2 (\bar P)$ is in the image of
$\iota$ and will not make distinction between $H^2(\pi_1 Y;\Z_2)$ and
its image in $H^2 (Y;\Z_2)$. Because of the identification of Proposition
\ref{P:hol}, the Casson invariant $\lambda'''(Y,w)$ can be defined by
counting points in the space $\P_c (Y;SU(2))$ with $[c] = w$, perhaps
after perturbation.


\section{The two--orbits}
According to the action of $H^1(Y;\Z_2) = (\Z_2)^3$ the space $\P_c(Y;SU(2))$
splits into orbits of possible orders one, two, four, and eight. In this
section we study the \emph{two--orbits} (orbits with two elements, or those
with stabilizer $\ZZ$).


\subsection{The two-orbits and invariant $\lambda'''$}

Consider a subgroup of $SO(3)$ that is isomorphic to $\ZZ$.  Such a
subgroup is generated by $180^{\circ}$ rotations about
two perpendicular axes in $\mathbb R^3$, and any two such subgroups
are conjugate to each other in $SO(3)$.  Hence the following
definition makes sense. Define $\R_w (Y;\ZZ)$ to be the subspace of
$\R_w (Y;SO(3))$ consisting of the $SO(3)$ conjugacy classes of
representations $\alpha: \pi_1(Y)\to SO(3)$ which factor through
$\ZZ\subset SO(3)$.

\begin{proposition}\label{P:two-orbits}
Let $[c] = w$ be a non--trivial class in $H^2(Y;\Z_2)$. Then the map $\pi:
\P_c(Y;SU(2))\to \R_w (Y;SO(3))$ establishes a bijective correspondence
between the set of two--orbits in $\P_c (Y;SU(2))$ and the set $\R_w(Y;\ZZ)$.
\end{proposition}

\begin{proof}
Suppose that the conjugacy class of a projective representation $\rho:
\pi_1(Y)\to SU(2)$ is fixed by a subgroup $\ZZ$ of $H^1(Y;\Z_2)$ generated
by homomorphisms $\alpha$, $\beta: \pi_1(Y) \to \Z_2$. Then there exists a
$u\in SU(2)$ such that $\alpha(x)\rho(x) = u\rho(x)u^{-1}$ for all $x\in
\pi_1(Y)$. Observe that $\rho(x) = u^2\rho(x)u^{-2}$ and, since $\rho$ is
irreducible, $u^2 = \pm 1$. The case $u^2 = 1$ should be excluded because
$u^2 = 1$ would imply that $u =\pm 1$ so that $-\rho(x) =\rho(x)$ at least
for some $x$, which is impossible in $SU(2)$. Therefore $u^2 = -1$ and,
after conjugation if necessary, we may assume that $u = i$. Then, for every
$x\in\pi_1(Y)$, we have $\pm\rho(x) = i\rho(x)i^{-1}$ so that $\im(\rho)
\subset S_i\cup j\cdot S_i$. Here, $S_i$ is the complex circle in $SU(2)$
(and $SU(2)$ is viewed as the group of unit quaternions).

Similarly, there exists a $v\in SU(2)$ such that $\beta(x)\rho(x) = v\rho(x)
v^{-1}$ and $v^2 = -1$. After conjugation by a complex number, we may assume
that $v = ia + bj$ where $a,b\in\mathbb R$ and $b\ge 0$. Next, $\alpha(x)
\beta(x)\rho(x) = (iv)\,\rho(x)\,(iv)^{-1}$ so that $(iv)^2 = -1$. An easy
calculation with quaternions shows that $v = j$ (and then $iv = k$). Thus
$\rho$ has the property that $\pm\rho(x) = i\,\rho(x)\,i^{-1}$ and $\rho(x)=
j\,\rho(x)\,j^{-1}$ for all $x\in\pi_1 (Y)$. Therefore,
\[
\im(\rho)\subset (S_i\cup j\cdot S_i)\,\cap\,(S_j\cup i\cdot S_j)
\]
where $S_j$ is the circle of quaternions of the form $\exp(j\phi)$. One can
easily see that the latter intersection is the group $Q = \{\,\pm 1, \pm i,
\pm j, \pm k\,\}$.

The above argument shows that any projective representation $\rho: \pi_1(Y)
\to SU(2)$ stabilized by $\ZZ$ factors through $Q$ and therefore its
associated $SO(3)$ representation $\ad\rho$ factors through a copy of
$\ZZ\subset SO(3)$.

To complete the proof, we only need to show that the orbit of $\rho$ consists
of exactly two points. Let $\gamma$ be a vector in $H^1(Y;\Z_2)$ completing
$\alpha$, $\beta$ to a basis. Then $\rho$ and $\rho^{\gamma}$ lie in the same
$H^1(Y;\Z_2)$--orbit but are not conjugate. The latter can be seen as follows:
if there exists a $w\in SU(2)$ such that $\gamma(x)\rho(x) = w\,\rho\,w^{-1}$
then $w = \pm k$ and $\alpha(x)\beta(x)\gamma(x)\rho(x) = (ijk)\,\rho(x)\,
(ijk)^{-1} = \rho(x)$ for all $x$, a contradiction.
\end{proof}

\begin{remark}\label{R:one-orbits}
The above proof shows in particular that no point of $\P_c (Y;SU(2))$ with
$[c]\ne 0$ is fixed by the entire group $H^1(Y;\Z_2)$ so that $\P_c(Y;SU(2))$
has no orbits of order one.
\end{remark}


\subsection{The number of two-orbits}

Our next goal is to find a formula for the number of points in
$\R_w (Y;\ZZ)$ modulo 2.

\begin{proposition}\label{P:cup}
Let $0\ne w\in H^2(Y;\Z_2)$ then
$\#\,\R_w (Y;\ZZ) = (a_1\cup a_2\cup a_3)\,[Y]\pmod 2$.
\end{proposition}

\begin{proof}
We begin by observing that any two subgroups of $SO(3)$ that are
isomorphic to $\ZZ$ are conjugate, and that moreover any automorphism of
such a subgroup is realized by conjugation by an element of $SO(3)$. Let
us fix a subgroup $\ZZ$ and a basis in it.

Since $\ZZ$ is abelian, every $\alpha\in\R_w (Y;\ZZ)$ factors through a
homomorphism $H_1(Y;\Z)\to \ZZ$. The two components of this homomorphism
determine elements $\beta,\gamma \in\Hom(H_1(Y);\Z_2)\cong H^1(Y;\Z_2)$.
It is straightforward to see that the $SO(3)$ representation $\alpha$ may
be recovered from $\beta$ and $\gamma$ via the formula $\alpha\cong \beta
\oplus \gamma \oplus \det(\beta \oplus \gamma)$. Since any element of
$\Lambda^2 H^1(Y;\Z_2)$ can be represented in the form $\beta\wedge\gamma$,
this establishes a one-to-one correspondence
\begin{equation}\label{E:lambda}
\R(Y;\ZZ) \to \Lambda^2 H^1(Y;\Z_2),
\end{equation}
where $\R(Y;\ZZ)$ is union of $\R_w (Y;\ZZ)$ over all possible $w$. Since
$H_1(Y;\Z)$ is torsion free, any element in $H^1(Y;\Z_2)$ is the mod 2
reduction of a class in $H^1(Y;\Z)$. It follows that the cup product of
any element $a\in H^1(Y;\Z_2)$ with itself is $0$. We compute
\begin{equation}\label{whitney}
\begin{split}
w_2(\alpha) & = w_1(\beta) \cup w_1(\gamma) +  w_1(\beta) \cup
w_1(\det(\beta \oplus \gamma)) + w_1(\gamma) \cup w_1 (\det(\beta
\oplus \gamma)) \\
& = w_1(\beta) \cup w_1(\gamma) + w_1(\beta) \cup (w_1(\beta) + w_1(\gamma)) +
w_1(\gamma) \cup (w_1(\beta) + w_1(\gamma))\\
&= w_1(\beta) \cup w_1(\gamma).
\end{split}
\end{equation}
Since $w_1(\beta) = \beta$ and $w_1(\gamma) = \gamma$, this shows that
$w_2(\alpha)$ is the image of $\beta\wedge\gamma$ under the map (\ref{E:cup}).

The result now follows by composing (\ref{E:cup}) and (\ref{E:lambda})\,:
if the triple cup product on $H^1(Y;\Z_2)$ vanishes mod 2 then the map
(\ref{E:cup}) is identically zero, hence $\R_w (Y;\ZZ)$ is empty for $w\ne
0$. If the triple cup product is nontrivial mod 2 then the map (\ref{E:cup})
is an isomorphism and $\R_w (Y;\ZZ)$ consists of exactly one element for
every choice of $0\ne w\in H^2(Y;\Z_2)$.
\end{proof}


\subsection{Non--degeneracy of the two--orbits}\label{S:4.3}

We wish to use Proposition \ref{P:cup} to calculate the contribution of the
two--orbits to $\lambda'''(Y)\pmod 2$. In order to do that, we need to
check the non--degeneracy condition for such orbits in the case when
the triple cup product on $Y$ is non--trivial.

\begin{proposition}\label{P:2nondegen}
Let $Y$ be an odd homology torus then $H^1 (Y;\ad\rho) = 0$ for any
projective representation $\rho: \pi_1(Y)\to SU(2)$ such that $\ad\rho
\in\R_w (Y;\ZZ)$ with $w\ne 0$.
\end{proposition}

To prove this proposition we notice that any $\ad\rho\in\R_w (Y;\ZZ)$
splits as $\ad\rho = \alpha_1\,\oplus\,\alpha_2\,\oplus\,\alpha_3$ where
each $\alpha_i: \pi_1(Y)\to\Z_2$ is a non--trivial representation, and
\[
H^1(Y;\ad\rho) = H^1(Y;\alpha_1)\,\oplus\,H^1(Y;\alpha_2)\,\oplus\,
H^1(Y;\alpha_3),
\]
compare with the proof of Proposition \ref{P:cup}. Let $\alpha: \pi_1(Y)
\to\Z_2 = \O(1)$ be a non--trivial representation, and $Y_{\alpha}$ the
regular double covering of $Y$ with $\pi_1(Y_{\alpha}) = \ker(\alpha)$.

\begin{lemma}\label{L:14}
The group $H^1(Y;\alpha) = H^1(Y;\mathbb R_{\,\alpha})$ is isomorphic to
the $(-1)$--eigen\-space of $\Z_2$ acting on $H^1(Y_{\alpha};\mathbb R)$.
\end{lemma}

\begin{proof}
This is immediate from the definition of $H^1(Y;\mathbb R_{\,\alpha})$.
\end{proof}

\begin{lemma}\label{L:five}
The cup product map $\cup\,a: H^1(Y;\Z_2)\to H^2(Y;\Z_2)$ is non-trivial
for some $a\in H^1(Y;\Z_2)$ if and only if $Y$ is an odd homology torus.
\end{lemma}

\begin{proof}
Suppose $a\cup b\neq 0$.  By Poincar\'e duality, there is $c\in H^1(Y;\Z_2)$
with $a\cup b\cup c\neq 0$.  Conversely, note that the cup product of any
three basis elements is the same as the cup product of any other three basis
elements. So if there is a non-zero cup product, extend $\{\,a\,\}$ to a
basis $\{\,a,b,c\,\}$ with $a\cup b\cup c\neq 0$. In particular $a\cup b\neq
0$ (and also $a \cup c\neq 0$.)
\end{proof}

\begin{remark}\label{R:five}
Note that we in fact proved that, if $Y$ is an odd homology torus, the rank
of $\cup\, x: H^1(Y;\Z_2)\to H^2(Y;\Z_2)$ equals two for any nonzero $x\in
H^1(Y;\Z_2)$.
\end{remark}

\begin{lemma}\label{L:15}
If $(a_1\cup a_2\cup a_3)\,[Y] = 1$ mod 2, then $H^1(Y_{\alpha};\Z) =
\Z^3$.
\end{lemma}

\begin{proof}
Let us consider the Gysin exact sequence for the double covering $\pi:
Y_{\alpha}\to Y$ (with coefficients in $\Z_2$)
\begin{multline}\notag
H^0(Y)\xrightarrow{\alpha} H^1(Y)\xrightarrow{\pi^*} H^1(Y_{\alpha})\\
\longrightarrow H^1(Y) \xrightarrow{\alpha} H^2(Y)\xrightarrow{\pi^*}
H^2(Y_{\alpha})\to H^2(Y)\xrightarrow{\alpha} H^3(Y)
\end{multline}
where the arrows marked $\alpha$ stand for the homomorphisms given by the
cup product with $w_1(\alpha)\in H^1(Y;\Z_2)$. This sequence works out to
\begin{multline}\notag
0\to\Z_2\xrightarrow{\alpha}(\Z_2)^3\xrightarrow{\pi^*} H^1(Y_{\alpha})\\
\longrightarrow (\Z_2)^3\xrightarrow{\alpha}(\Z_2)^3\xrightarrow{\pi^*}
H^2(Y_{\alpha})\to(\Z_2)^3 \xrightarrow{\alpha} \Z_2 \to 0
\end{multline}
The first and the last zeroes are because $w_1(\alpha)\ne 0$. According to
Remark \ref{R:five}, the image of $\cup w_1(\alpha): (\Z_2)^3\to (\Z_2)^3$
has rank two. Counting ranks we get that $H^1(Y_{\alpha};\Z_2) =
H^2(Y_{\alpha};\Z_2) = (\Z_2)^3$. The result now follows by the universal
coefficient theorem.
\end{proof}

\begin{proof}[Proof of Proposition~\ref{P:2nondegen}] Let us fix an 
isomorphism between integral homology of $Y$ and that of the
$3$-torus $T$, and choose a map $f: Y\to T$ that induces this
isomorphism. Let $T_{\alpha}$ be the double covering of $T$ corresponding
to the homomorphism $\alpha: H_1 (T;\Z)\to\Z_2$ which makes the following
diagram commute

\[
\begin{CD}
\pi_1 (Y) @>\alpha >> \Z_2 \\
@VVV @| \\
H_1(T;\Z) @>\alpha >> \Z_2
\end{CD}
\]

\bigskip\noindent
The map $\pi_1(Y)\to H_1(T;\Z)$ in this diagram is obtained by composing the
abelianization $\pi_1(Y)\to H_1(Y;\Z)$ with the isomorphism $H_1(Y;\Z) =
H_1(T;\Z)$ induced by $f$.

Let $f_{\alpha}: Y_{\alpha}\to T_{\alpha}$ be a lift of $f$. Comparing Gysin
exact sequences of $\pi: Y_{\alpha}\to Y$ and $\pi: T_{\alpha}\to T$ using
the five--lemma, we conclude that the map $f_{\alpha}^*: H^1(T_{\alpha};\Z)
\to H^1(Y_{\alpha};\Z)$ is an isomorphism when tensored with $\Z_2$. Since
$H^1(Y_{\alpha};\Z) = \Z^3$ by Lemma \ref{L:15}, we also conclude that
$f_{\alpha}^*: H^1(T_{\alpha};\mathbb R)\to H^1(Y_{\alpha};\mathbb R)$ is an
isomorphism.

This implies that $\pi^*: H^1(Y;\mathbb R)\to H^1(Y_{\alpha};\mathbb R)$ is an
isomorphism, for this is true for $Y= T$, and we just observed the isomorphism
in the upper line of the following commutative diagram

\[
\begin{CD}
H^1(T_{\alpha};\mathbb R) @>= >> H^1(Y_{\alpha};\mathbb R) \\
@A\pi^* A= A  @A\pi^* AA \\
H^1(T;\mathbb R) @>= >> H^1(Y;\mathbb R)
\end{CD}
\]

\bigskip\noindent
On the other hand, the image of $\pi^*$ equals the $(+1)$--eigenspace of the
$\Z_2$--action on $H^1(Y_{\alpha};\mathbb R)$. Together with Lemma \ref{L:14},
this implies that $H^1(Y;\alpha) = H^1(Y;\mathbb R_{\alpha}) = 0$.
\end{proof}


\section{Perturbations}\label{S:perturb}

In this section we deal with the situation when the critical point set
$\M(P)$ of $\cs: \B(P) \to \mathbb R/\Z$ is degenerate. We describe a
class of equivariant admissible perturbations $h: \B(P) \to \mathbb R$
and prove that, for a small generic $h$, the critical point set of $\cs
+h$ is non-degenerate. This set is acted upon by $H^1 (Y;\Z_2)$ in such
a manner that an argument similar to that used in the non--degenerate
case completes the proof of Theorem \ref{T:main}. Following the approach
originated by Taubes, Floer and
Donaldson~\cite{taubes:casson,floer:instanton,donaldson:floer}, we
use holonomy around loops to define a perturbation of the function
$\cs$.  The main issue is to choose loops so that the perturbation is
$H^1 (Y;\Z_2)$-equivariant; this is done by imposing a simple
homological restriction on the loops.  We largely employ the notation
of C.~Herald~\cite{herald:perturbations}.


\subsection{Equivariant admissible perturbations}
Let $\gamma_k: S^1\times D^2 \to Y$, $k = 1,\ldots, n$, be a collection
of embeddings of solid tori in $Y$ with disjoint images. We will use
the same notation for the loops $\gamma_k (S^1\times\{\,0\,\})$ in $Y$,
and call $\Gamma = \{\,\gamma_k\,\}$ a \emph{link}. A link $\Gamma$ is
called \emph{mod--$2$ trivial} if $0 = [\gamma_k] \in H_1 (Y;\Z_2)$ for
all $k$. Let $\eta(z)$ be a smooth rotationally symmetric bump function
on the unit disc $D^2$ with support away from the boundary of $D^2$
and with integral one. Finally, let $f_k: SU(2)\to \mathbb R$, $k = 1,
\ldots,n$, be $C^3$--functions invariant with respect to conjugation.

Following the construction in Section~\ref{S:hol}, choose a lifting
of the holonomy to $SU(2)$. It is uniquely determined by a choice of
square roots of $\hol_C$ on a set of representative loops, different
choices leading to equivalent theories. For each based loop $\gamma$,
we obtain a well defined map $\A(P)\to SU(2)$. Define
\begin{equation}\label{ha}
h(A) = \sum_{k=1}^n\;\int_{D^2}\;f_k(\hol_A(\gamma_k(S^1\times\{\,z\,\})))
\,\eta(z)\,d^2 z,
\end{equation}
where $\hol_A (\gamma_k(S^1\times\{\,z\,\}))$ stands for holonomy of $A$
around the loop $\gamma_k (S^1\times\{\,z\,\})$, $z\in D^2$. A basepoint
must be chosen in order for this holonomy to be well defined; however,
the conjugation invariance of $f_k$ ensures that the function $h$ does
not depend on such choices. The action of $\G(P)$ only changes holonomies
around $\gamma_k (S^1\times\{\,z\,\})$ within their $SU(2)$ conjugacy
class. Thus we have a well defined function
\begin{equation}\label{hb}
h: \B(P)\to\mathbb R
\end{equation}
which we call an \emph{admissible perturbation relative to $\Gamma$}.
For any link $\Gamma$, denote by $\H_{\Gamma}$ the space of admissible
perturbations relative to $\Gamma$ with the $C^3$--topology given by
the correspondence $h \mapsto (f_1,\ldots, f_n)$.

\begin{lemma}\label{lemma5.1}
If $\Gamma$ is mod--$2$ trivial then the function $h$ defined in
\eqref{hb} is $H^1 (Y;\Z_2)$--invariant.
\end{lemma}

\begin{proof}
We need to prove that $h(A\otimes\chi) = h(A)$ for all $A\in\A(P)$
and $\chi\in H^1(Y;\Z_2)$. This follows easily from the formula
\[
\hol_{A\otimes\chi}(\gamma_k(S^1\times\{\,z\,\})) =
\hol_A(\gamma_k(S^1\times\{\,z\,\}))\cdot\chi(\gamma_k)
\]
after we notice that $\chi(\gamma_k) = 1$ because $\chi: \pi_1(Y)
\to \Z_2$ factors through $H_1(Y;\Z_2)$.
\end{proof}

Any admissible perturbation $h \in \H_{\Gamma}$ where $\Gamma$ is a
mod--$2$ trivial link will be called an \emph{equivariant admissible
perturbation}.


\subsection{Perturbed projectively flat connections}

Let $h: \A(P)\to\mathbb R$ be an admissible perturbation relative to a
link $\Gamma$. The projection map $\pi: \A(P)\to\A(\bar P)$ identifies
the tangent space of $\A(P)$ with that of $\A(\bar P)$. Identify the
latter with $\Omega^1 (Y; \ad\bar P)$ and define $\zeta_h: \A(P) \to
\Omega^1 (Y;\ad\bar P)$ by the formula
\[
\zeta_h (A) = *\,F_{\bar A} - 4\pi^2\cdot \nabla h (A),
\]
where $\nabla h$ is the $L^2$--gradient of $h$. A straightforward
calculation shows that, up to the identification of the tangent spaces,
$\zeta_h$ is just $-4\pi^2$ times the $L^2$--gradient of the function
$\cs + h$.

A connection $A\in\A(P)$ is called \emph{$h$--perturbed projectively
flat} if $\zeta_h (A) = 0$. The moduli space of $h$--perturbed
projectively flat connections is denoted by $\M_h(P)$ so that $\M_h(P)
= \zeta_h^{-1}(0)/\G(P)$. If $h = 0$ then $\M_h (P)$ coincides with
the moduli space $\M(P)$ of projectively flat connections, see Section
\ref{S:2.2}.

Next we wish to describe the local structure of $\M_h(P)$ near a point
$[A]\in\M_h(P)$. The slice through $A$ to the $\G(P)$--action on $\A(P)$
is the affine subspace
\[
X_A = \{\,A + \pi_*^{-1}(a)\;|\; a\in\ker d_{\bar A}^*\,\}\subset\A(P)
\]
where $d_{\bar A}^*: \Omega^1 (Y;\ad\bar P) \to \Omega^0 (Y;\ad\bar P)$.
Since $c_1(P)$ is an odd element in $H^2 (Y;\Z)$, the stabilizer of $A$
in $\G(P)$ coincides with the center of $SU(2)$, and a small
neighborhood of $A$ in $X_A$ gives a local model for $\B(P)$ near $[A]$.
Therefore, the moduli space $\M_h (P)$ near $[A]\in\M_h (P)$ is the zero
set of $\zeta_h$ restricted to the slice $X_A$. The linearization of
$\zeta_h: \A(P)\to\Omega^1(Y;\ad\bar P)$ at $A\in\A(P)$ is the operator
\[
*d_{A,h} = *d_{\bar A} - 4\pi^2\cdot\Hess h (A): \Omega^1 (Y;\ad\bar P)
\to \Omega^1 (Y;\ad\bar P),
\]
hence the tangent space to $\M_h (P)$ at $A \in \M_h (P)$ can be
identified with
\[
H^1_h (Y;\ad A) = \ker *d_{A,h}/\,\im\{\,d_{\bar A}: \Omega^0 (Y;
\ad\bar P) \to \Omega^1 (Y; \ad\bar P)\,\}.
\]

We call $\M_h (P)$ \emph{non--degenerate} at $[A] \in \M_h(P)$ if
$H^1_h(Y;\ad A) = 0$; we call it \emph{non--degenerate} if it
is non--degenerate at all $[A] \in \M_h (P)$. If $\M_h(P)$ is
non--degenerate, it consists of finitely many points, and their
signed count gives $\lambda'''(Y,w)$ where $c_1(P) = w\pmod 2$.

If $h$ is an equivariant admissible perturbation then according to 
Lemma~\ref{lemma5.1}, $\M_h (P)$ is acted
upon by $H^1 (Y;\Z_2)$.


\subsection{Abundance of equivariant admissible perturbations}
Our main goal in the next few sections will be to show that one can
always find an equivariant admissible perturbation $h$ such that $\M_h(P)$
is non--degenerate. We begin by choosing a mod--$2$ trivial link $\Gamma$
satisfying certain necessary conditions. Such links are called abundant;
the definition of abundance at $A$ depends on the size of the stabilizer
of $A$ in $H^1 (Y;\Z_2)$.

Let $\Gamma = \{\,\gamma_k\,\}$ be a mod--$2$ trivial link and $A$ a
projectively flat connection whose stabilizer in $H^1(Y;\Z_2)$ is trivial.
Then $\Gamma$ is called \emph{abundant} at $A$ if there exist admissible
perturbations $h_i \in \H_{\Gamma}$, $i = 1,\ldots, m$, such that the map
from $\mathbb R^m$ to $\Hom(H^1(Y;\ad A),\mathbb R)$ given by
\begin{equation}\label{E:holmap}
(x_1,\ldots,x_m)\mapsto\sum\; x_i\, Dh_i (A)
\end{equation}
is surjective.

Now, let $A$ be a projectively flat connection whose stabilizer in
$H^1(Y;\Z_2)$ equals $\Z_2$. Let $\tau$ be a generator in $\Z_2$ and
denote by $V^{\pm}(A)$ respectively the $(\pm 1)$--eigenspaces of
$\tau_*: H^1 (Y;\ad A)\to H^1 (Y;\ad A)$. Denote by $\Sym(V)$
the set of symmetric bilinear forms on a vector space $V$. A mod--$2$
trivial link $\Gamma$ is called \emph{abundant} at $A$ if there exist
admissible perturbations $h_i \in \H_{\Gamma}$, $i = 1,\ldots, m$,
such that the map from $\mathbb R^m$ to $\Hom(V^+(A),\mathbb R)\,
\oplus\,\Sym(V^-(A))$ given by

\[
(x_1,\ldots,x_m)\mapsto \left(\sum\; x_i\,Dh_i (A),\;\sum\; x_i\,
\Hess h_i(A)\right)
\]
is surjective (this definition makes sense because $h$ is
$H^1(Y;\Z_2)$--equivariant).

Due to the fact (cf. Section \ref{S:4.3}) that $H^1(Y;\ad A) = 0$ for
any projectively flat connection $A$ whose stabilizer is bigger than
$\Z_2$, we do not need to perturb $A$ and hence do not need the concept
of abundancy at such a connection.

A useful remark is that if $\Gamma $ is abundant at $A$, and $\Gamma_0$
is another link whose components are close to those of $\Gamma$, then
$\Gamma_0$ is also abundant at $A$. Moreover, the perturbation functions
$h_i$ can be taken to be the same as for $\Gamma$. These facts come from
the homotopy invariance of parallel transport. Note also that if $\Gamma$
is abundant at $A$ and $\Gamma \subset \Gamma'$ then $\Gamma'$ is also
abundant at $A$. The following result will be proved in Section
\ref{S:lemmas} below.

\begin{proposition}\label{P:ab}
There exists a mod--$2$ trivial link $\Gamma$ which is abundant at all
$A \in \M(P)$ whose stabilizer is at most $\Z_2$.
\end{proposition}


\subsection{Non--degeneracy results}
In this section, we will make use of Proposition \ref{P:ab} to prove
existence of equivariant admissible perturbation functions making $\M_h(P)$
non--degenerate.

Let $\Gamma$ be an abundant mod--$2$ trivial link as in Proposition
\ref{P:ab}. Let $\B^*$ be the subset of $\B(P)$ consisting of connections
whose stabilizer in $H^1 (Y;\Z_2)$ is trivial, and consider the universal
zero set
\[
\mathcal Z^* = \{\,([A],h)\in\B^*\times\H_{\Gamma}\;|\;\zeta_h(A) = 0\,\}.
\]
The moduli space $\M^* = \M(P) \cap \B^*$ will be viewed as a subset of
$\mathcal Z^*$ by assigning $([A],0)$ to $[A] \in \M^*$. The following
proposition roughly states that $\M^*$ can be ``thickened'' inside
$\mathcal Z^*$ to become a manifold.

\begin{proposition}\label{P:star}
The moduli space $\M^*$ has an open neighborhood $\mathcal U^*\subset
\mathcal Z^*$ which is a submanifold of $\B^*\times\H_{\Gamma}$.
\end{proposition}

\begin{proof}
Fix a point $[A_0] \in \M^*$ and consider the map
\[
P: X_{A_0}\times\H_{\Gamma}\to \ker d^*_{\bar A_0}
\]
given by $P(A,h) = \Pi_{A_0}\zeta_h (A)$ where $\Pi_{A_0}: \Omega^1
(Y;\ad\bar P) \to \ker d^*_{\bar A_0}$ is the $L^2$--orthogonal
projection. The first partial derivative of this map is Fredholm with
cokernel $H^1(Y;\ad A_0)$. Since $\Gamma$ is abundant at $A_0$,
the image of the partial derivative $\p P/\p h$ is a subspace which
orthogonally projects onto this cokernel. Therefore, $P$ is a
submersion at $[A_0]$. The implicit function theorem now implies
that $P^{-1}(0)\subset X_{A_0}\times\H_{\Gamma}$ is smooth near $(A_0,
0)$. Moreover, $\Pi_{A_0}\zeta_h(A) = 0$ if and only if $\zeta_h (A)
= 0$, at least in a small neighborhood of $A_0$ in $X_{A_0}$, see
\cite[Lemma 12.1.2]{morgan:mrowka:ruberman}. The union of such
neighborhoods over all $[A_0] \in \M^*$ is the open submanifold
$\mathcal U^*$.
\end{proof}

\begin{corollary}\label{C:one}
For a small generic perturbation $h\in\H_{\Gamma}$, the moduli space
$\M_h^* = \M_h(P)\cap \B^*$ is non-degenerate.
\end{corollary}

\begin{proof}
The projection from $\mathcal U^*$ to $\H_{\Gamma}$ is Fredholm of
index zero hence the result follows from the Sard--Smale theorem.
\end{proof}

Let us now turn to connections in $\B(P)$ with stabilizer $\Z_2$. Fix
a generator $\tau$ in a copy of $\Z_2$ and consider the subset
$\B^{\tau}$ of $\B(P)$ consisting of gauge equivalence classes of
connections stabilized by $\tau$. The argument of Proposition
\ref{P:star}, after a slight modification,  can be used to
prove that $\M_h^{\tau} = \M_h(P) \cap \B^{\tau}$ is non-degenerate
\emph{inside $\B^{\tau}$} for a generic small perturbation $h \in
\H_{\Gamma}$. However, we are interested in non-degeneracy inside $\B(P)$
and hence in a description of the normal bundle of $\M_h^{\tau}$ in
$\M_h(P)$.

To describe this normal bundle, we need to review the Kuranishi model of
$\M(P)$ near $[A]\in \M^{\tau}$, see \cite{herald:perturbations}. Since
the derivative of $\Pi_A \zeta(A): X_A\to \ker d^*_{\bar A}$ is already
a Fredholm isomorphism from the orthogonal complement of its kernel to
the orthogonal complement of its cokernel, the effect on the normal bundle
of adding a small perturbation $h$ is determined by $\Hess h (A): V^- (A)
\to V^- (A)$. In particular, the normal bundle is zero dimensional whenever
$\Hess h (A)$ is an isomorphism.

Let us consider the universal zero set
\[
\mathcal Z^{\tau} = \{\,([A],h)\in\B^{\tau}\times\H_{\Gamma}\;|\;\zeta_h(A)
= 0\,\}
\]
and view $\M^{\tau} = \M(P)\,\cap\,\B^{\tau}$ as a subset of
$\mathcal Z^{\tau}$ by assigning $([A],0)$ to every $[A] \in \M^{\tau}$.

\begin{proposition}
The moduli space $\M^{\tau}$ has an open neighborhood $\mathcal U^{\tau}
\subset \mathcal Z^{\tau}$ such that
\begin{enumerate}
\item[(a)]
$\mathcal U^{\tau}$ is a submanifold in $\B^{\tau}\times\H_{\Gamma}$, and
\item[(b)]
for every $A$, a generic $h$ such that $([A],h)\in\mathcal U^{\tau}$, has
non-degenerate Hessian.
\end{enumerate}
\end{proposition}

\begin{proof}
Let us fix $[A_0] \in \M^{\tau}$. The slice at $A_0$ of the gauge group
action on $\B^{\tau}$ is given by
\[
X^{\tau}_{A_0} = \{\,A_0+\pi_*^{-1}(a)\;|\;a\in\ker d^*_{\bar{A_0}}\cap
\Omega^1 (Y;\ad\bar P)^+\,\},
\]
where $\Omega^1 (Y;\ad\bar P)^+$ is the $(+1)$--eigenspace of
$\tau: \Omega^1(Y; \ad\bar P) \to \Omega^1(Y; \ad\bar P)$. Denote by
$\underline{\Sym }(V^-)$ the bundle over an open neighborhood $W$ of
$(A_0,0)$ in $X^{\tau}_{A_0}\times\H_{\Gamma}$ whose fiber over
$(A,h)$ is $\Sym(V_h^-(A))$, the set of symmetric bilinear forms on
the $(-1)$--eigenspace $V_h^-(A)$ of $\tau: H^1_h (Y; \ad A)\to
H^1_h (Y; \ad A)$). Let
\[
P: W\to (\ker d^*_{\bar{A_0}}\cap\Omega^1 (Y;\ad\bar P)^+)\,\oplus\,
\underline{\Sym } (V^-)
\]
be the section $P(A,h) = (\Pi^*_{A_0}\zeta_h(A),\Hess h (A))$ where
$\Pi^*_{A_0}$ is $\Pi_{A_0}$ followed by the $L^2$--orthogonal
projection onto $\ker d^*_{\bar{A_0}}\cap\Omega^1(Y;\ad\bar P)^+$.
The first partial derivative of $P$ at $(A_0,0)$ has cokernel
$V^+ (A_0)\,\oplus\,\Sym(V^-(A_0))$. Since $\Gamma$ is
abundant at $A_0$, the image of the partial derivative $\p P/
\p h$ is a subspace which orthogonally projects onto this cokernel.
The implicit function theorem now implies that $P^{-1}(\{\,0\,\}
\times \Sym(V^-(A_0)))$ is smooth near $(A_0,0)$, which proves
part (a). Since non-degenerate symmetric forms are generic in
$\Sym(V^-(A_0))$, the part (b) also follows.
\end{proof}

\begin{corollary}\label{C:two}
For a small generic admissible perturbation $h\in\H_{\Gamma}$ the moduli
space $\M^{\tau}_h = \M_h(P)\,\cap\,\B^{\tau}$ is non--degenerate.
\end{corollary}

\begin{proof}
The projection from $\mathcal U^{\tau}$ to $\H_{\Gamma}$ is Fredholm of
index zero hence the result follows from the Sard--Smale theorem.
\end{proof}

\subsection{Proof of Proposition \ref{P:ab}}\label{S:lemmas}
The proof of Proposition \ref{P:ab} naturally divides into three
parts, which can be viewed as pointwise, local, and global abundance.
The passage from pointwise to local and global abundance is proved in
essentially the same manner  as in~\cite{herald:perturbations}.
These rely on basic analytical properties of the Chern-Simons
function, especially the compactness of the perturbed moduli space,
and the restriction to equivariant admissible perturbations does not
change these arguments.  Thus we will concentrate on establishing
pointwise abundance, as summarized in the following proposition.

\begin{proposition}\label{lm1}
If $A$ is a projectively flat connection then there exists a mod--$2$
trivial link $\Gamma$ which is abundant at $A$. Moreover, if the
stabilizer of $A$ equals $\Z_2$, the functions $h_1,\ldots, h_m$ can
be chosen so that, for some $k$,
\begin{enumerate}
\item[(a)] $Dh_1(A),\ldots, Dh_k(A)$ span $\Hom(V^+,\mathbb R)$,
\item[(b)] $\Hess h_{k+1}(A),\ldots, \Hess h_m(A)$ span $\Sym (V^-)$,
and
\item[(c)] $Dh_j(A) = 0$ for $j = k+1,\ldots,m$.
\end{enumerate}
\end{proposition}
For connections with trivial stabilizer, the result is established in
Lemma~\ref{8orbit}, while the result for connections with stabilizer
$\Z_2$ is in Lemma~\ref{L3:4orbit}.

Denote by $p: \yhom \to Y$ the regular covering space corresponding
to the surjection $\phi_2: \pi_1(Y) \to H_1(Y;\Z_2) \cong (\Z_2)^3$.
This cover might be called the $2$--universal abelian cover, because
of the following observation.  Let $G$ be an abelian group which is
a $\Z_2$--vector space, and suppose that $f: \pi_1(Y) \to G$ is a
homomorphism.  Then there is a unique homomorphism $\hat f:
H_1(Y;\Z_2) \to G$ such that $\hat f \circ \phi_2 = f$.  This can be
readily seen from the universal property of the abelianization
$\phi:\pi_1(Y) \to H_1(Y;\Z)$ and the universal property of the map
$H_1(Y;\Z) \to H_1(Y;\Z) \otimes \Z_2 \cong H_1(Y;\Z_2)$.

For a connection $A$ on the bundle $P \to Y$, we will denote by
$\ahom$ its pull-back to $\yhom$. We need to understand the
behavior of a projectively flat connection on $Y$, when lifted in
this manner to $\yhom$.

\begin{lemma}\label{pullback}
Let $\rho: \pi_1 Y\to SU(2)$ be a projective representation and
$\rhom: \pi_1\yhom \to SU(2)$ the induced projective representation.
Let $\Stab(\rho)$ denote the stabilizer of $\rho$ in $H^1(Y;\Z_2)$
then
\begin{enumerate}
\item[(a)] $\Stab(\rho) = 1$ if and only if $\rhom$ is irreducible,
\item[(b)] $\Stab(\rho) = \Z_2$ if and only if $\rhom$ is reducible
non-central, and
\item[(c)] $\Stab(\rho) = \ZZ$ if and only if $\rhom$ is central.
\end{enumerate}
No other stabilizers $\Stab(\rho)$ may occur.
\end{lemma}

\begin{proof}
According to Remark \ref{R:one-orbits}, the only $\Stab(\rho)$ that
occur are 1, $\Z_2$, and $\ZZ$.

Suppose that $\Stab(\rho) = \ZZ$ then, as we saw in the proof of
Proposition \ref{P:two-orbits}, the image of $\rho$ is contained
in a copy of the group $Q = \{\,\pm 1, \pm i, \pm j,\pm k\,\}$.
Since $\pi_1 \yhom$ is in the kernel of the map $\phi_2$, we
conclude that the image of $\rhom$ is contained in the kernel of
the corresponding map $Q \to H_1(Q;\Z_2) \cong \ZZ$. This kernel
is the same as the commutator subgroup $[Q,Q] = \{\,\pm 1\,\}$
hence $\rhom$ is central. Conversely, if $\rhom$ is central, its
adjoint representation $\ad\rhom$ is trivial so that $\im(\ad\rho)$
is contained in a subgroup of $SO(3)$ of order at most eight.
Therefore, $\im(\ad\rho)$ is contained in a copy of $\ZZ$, and
then $\im(\rho)\subset Q$. In particular, $\Stab(\rho) = \ZZ$.

Now suppose that $\Stab(\rho) = \Z_2$. As in the proof of Proposition
\ref{P:two-orbits} we see that the image of $\rho$ is contained in a
copy of $S_i\cup j\cdot S_i$ where $S_i$ is the unit complex circle.
By the argument about the abelianization mod $2$, it follows that
$\im(\rhom)\subset S_i$, so that $\rhom$ is abelian. Conversely, if
$\rhom$ is abelian then $\im(\ad\rhom)$ is contained in a copy of
$SO(2)$, and $\im(\ad\rho)$ in its finite 2--prime extension.
Therefore, $\im(\rho)$ is contained in a copy of $S_i\cup j\cdot S_i$.

The remaining case follows by elimination.
\end{proof}
%
%
%

The same result holds for projectively flat connections in place of
projective representations.

Now we are able to deduce the existence of abundant links in the
simplest case, when the stabilizer of $A$ in $H^1(Y;\Z_2)$ is
trivial.
\begin{lemma}\label{8orbit}
Let  $A$ be a projectively flat connection whose stabilizer in
$H^1(Y;\Z_2)$ is trivial.  Then there exists a mod--$2$ trivial link
$\Gamma$ that is abundant at $A$.
\end{lemma}

\begin{proof}
By Lemma~\ref{pullback}, the connection $\ahom$ is irreducible,
and so by~\cite[Lemma 60]{herald:perturbations} there is a link
$\tilde\Gamma$ in $\yhom$ that is abundant at $\ahom$.  If we
perturb $\tilde\Gamma$ by a small amount, its projection $\Gamma =
p(\tilde\Gamma)$ will be an embedded link in $Y$.  It is clear
that $\Gamma$ is mod--$2$ trivial; we claim that in fact it is
abundant.  In the discussion that follows, the perturbing functions
on $Y$ will be the push-down of the perturbing functions $h_i$ on
$\yhom$.  This makes sense because the holonomy of $\ahom$ around
a component $\tilde\gamma$ of $\tilde\Gamma$ is the same as the
holonomy of $A$ around $p(\tilde\gamma)$.

Consider the commutative diagram

\[
\begin{CD}
\mathbb R^m @>>> \Hom(H^1(\yhom;\ad\ahom);\mathbb{R})\\
@| @VV{(p^*)^*}V\\
\mathbb R^m @>>> \Hom(H^1(Y;\ad A);\mathbb{R})
\end{CD}
\]

\medskip\noindent
where the horizontal arrows are the holonomy maps as in~\eqref{E:holmap}.
The arrow along the top is surjective, because $\tilde\Gamma$ is abundant
at $\ahom$.   Now it is a standard consequence of the transfer
sequence~\cite{brown} that the map
\begin{equation}\label{E:transfer8}
p^*: H^1(Y;\ad A) \to H^1(\yhom;\ad\ahom)
\end{equation}
is injective. Hence the bottom arrow is surjective as well.
\end{proof}

We next turn our attention to the abundance at projectively flat
connections with stabilizer $\Z_2$ in $H^1 (Y;\Z_2)$. Let $\rho:\pi_1 Y
\to SU(2)$ be a projective representation with $\Stab(\rho) = \Z_2$ and
fix a generator $\tau \in \Z_2 \subset H^1 (Y; \Z_2)$. Then $\tau$ acts
on $\P_c (Y; SU(2))$ fixing $\rho$ and hence inducing a $\Z_2$--action
$\tau^*$ on the tangent space $T_{\rho}\P_c(Y;SU(2)) = H^1 (Y;\ad\rho)$.
Denote as before by $V^{\pm}(\rho)$ the $(\pm 1)$--eigenspaces of
$\tau^*$ so that $H^1 (Y;\ad\rho) = V^+ (\rho)\,\oplus\, V^-(\rho)$.

According to Lemma \ref{pullback} the lift $\rhom: \pi_1\yhom\to SU(2)$
of $\rho$ is a reducible (non--central) projective representation.
Assuming (after conjugation if necessary) that $\im(\rhom)$ is contained
in the complex circle $S_i$, we obtain a splitting $\ad\rhom = \mathbb R
\,\oplus\, \ad_{\mathbb C}\rhom$ where $\mathbb R$ stands for a
trivial
one-dimensional representation and $\ad_{\mathbb C}\rhom: \pi_1 \yhom
\to SO(2)$. Accordingly, $H^1 (\yhom; \ad\rhom)$ splits as
\[
H^1 (\yhom; \ad\rhom) = H^1 (\yhom; \mathbb R) \,\oplus\,
H^1 (\yhom; \ad_{\mathbb C}\rhom).
\]

\begin{lemma}\label{L1:4orbit}
The projection $p: \yhom\to Y$ induces a monomorphism $p^*: H^1 (Y;\ad\rho)
\allowbreak \to H^1 (\yhom; \ad\rhom)$ such that
\[
p^* (V^+(\rho)) \subset H^1 (\yhom;\mathbb R)\quad\text{and}\quad
p^* (V^-(\rho))\subset H^1 (\yhom;\ad_{\mathbb C}\rhom).
\]
\end{lemma}

\begin{proof}
That $p^*$ is a monomorphism follows from the standard transfer argument,
see \cite{brown}. Since the conjugacy class of $\rho$ is fixed by $\tau$,
there exists an element $u\in SU(2)$ such that $u^2 = -1$ and $\tau(x)
\rho(x) = u\rho(x)u^{-1}$ for all $x\in\pi_1 Y$. If $x\in\pi_1\yhom$ then
$\rhom(x)= \rho(x)$ and $\tau(x)= 1$ so that $\rhom(x)= u\rhom(x)u^{-1}$.
Since $\rhom(x)\in S_i$ we conclude that $u = \pm i$.

To describe the induced action $\tau^*$ on $T_{\rho} \P_c (Y;SU(2)) =
H^1 (Y; \ad\rho)$ we first identify the tangent spaces at $\rho$ and
$\rho^{\tau}$ by $\ad u$, and then linearize the map $\rho \mapsto
\rho^{\tau}$ as follows\,:
\begin{multline}\notag
(1 + \ep\cdot\xi(x))\rho(x) \xrightarrow{\tau} \\
\tau(x)(1 + \ep\cdot\xi(x))\rho(x) = (1 + \ep\cdot\xi(x))\tau(x)\rho(x)
\xrightarrow{\ad u}  \\
u(1 + \ep\cdot\xi(x))\tau(x)\rho(x)u^{-1} =
(1 + \ep\cdot u\,\xi(x) u^{-1})\rho(x).
\end{multline}
Here, $\xi: \pi_1 Y\to \su(2)$ is a 1-cocycle representing an element of
$H^1 (Y; \ad\rho)$, and $\ep$ is a small positive real number. Thus the
action $\tau^*: H^1 (Y;\ad\rho) \to H^1 (Y;\ad\rho)$ at the level of
1-cocycles is given by the formula $\tau^*(\xi) = u \xi u^{-1}$. Since
$u=\pm i$, the subspace $V^+(\rho)$ is generated by 1-cocycles $\xi$ with
$\im(\xi) \subset i\,\mathbb R$, and $V^-(\rho)$ by 1-cocycles $\xi$ with
$\im(\xi)$ in the subspace $\mathbb C \subset \su(2)$ spanned by $j$ and
$k$.

The embedding $p^*: H^1 (Y; \ad\rho) \to H^1 (\yhom; \ad\rhom)$ is
given by pulling back the 1-cocycles $\xi: \pi_1 Y\to \su(2)$ via the
homomorphism
$p_*: \pi_1 \yhom\to \pi_1 Y$. In particular, if $\im(\xi)\subset i\,\mathbb
R$ then $\im(p^*\xi)\subset i\,\mathbb R$ so that $[p^*\xi]\in H^1 (\yhom;
\mathbb R)$. Similarly, if $\im(\xi)$ belongs to $\mathbb C$ spanned by $j$
and $k$ then $[\im(p^*\xi)]\in H^1(\yhom;\ad_{\mathbb C}\rhom)$.
\end{proof}

Denote by $\tilde\tau:\yhom\to\yhom$ the covering transformation
corresponding to the (dual of) $\tau\in\Z_2\subset H^1 (Y;\Z_2)$,
and by $\tilde\tau^*: H^1(\yhom;\ad\rhom)\to H^1 (\yhom;\ad\rhom)$
the induced action.

\begin{lemma}\label{L2:4orbit}
The subset $p^*(V^+(\rho))\subset H^1(\yhom;\mathbb R)$ is the
$(+1)$--eigenspace of $\tilde\tau^*: H^1 (\yhom;\mathbb R)\to H^1 (\yhom;
\mathbb R)$, and $p^*(V^-(\rho))\subset H^1 (\yhom;\ad_{\mathbb C}\rhom)$
the $(+1)$--eigenspace of $\tilde\tau^*: H^1 (\yhom;\ad_{\mathbb C}\rhom)
\to H^1 (\yhom;\ad_{\mathbb C}\rhom)$. Moreover, $p^*(V^-(\rho))$ is a
totally real subspace of the complex vector space $H^1 (\yhom;
\ad_{\mathbb C}\rhom)$.
\end{lemma}

\begin{proof}
The first two statements follow from the standard transfer argument,
see \cite{brown}. For the last statement, note that the pullback of
$\ad_{\mathbb C}\rhom$ via $\tilde\tau$ is exactly the complex
conjugate representation $\overline{\ad_{\mathbb C}\rhom}$.  It
follows that the action of $\tilde\tau^*$ on the cochains used to
compute $H^1 (\yhom;\ad_{\mathbb C}\rhom)$ is complex anti-linear,
and so the action on this cohomology group is also complex
anti-linear.  Thus the $(+1)$--eigenspace $p^*(V^-(\rho))$ is
totally real.
\end{proof}

\begin{lemma}\label{L3:4orbit}
Let  $A$ be a projectively flat connection whose stabilizer in
$H^1(Y;\Z_2)$ is $\Z_2$. Then there exists a mod--$2$ trivial link
$\Gamma$ that is abundant at $A$.
\end{lemma}

\begin{proof}
By Lemma~\ref{pullback}, the pull back connection $\ahom$ and the
associated projective representation $\rhom$ are reducible and
non--central. By \cite[Corollary 64 and Corollary 66]{herald:perturbations},
see also \cite[Proposition 3.4]{boden:herald:su3}, there is a link
$\tilde\Gamma$ in $\yhom$ and admissible perturbations $\tilde h_i:
\B(\tilde P)\to\mathbb R$, $i = 1,\ldots, m$, such that the map
\[
\mathbb R^m \to \Hom(H^1(\yhom;\mathbb R),\mathbb R)\,\oplus\,
\Herm (H^1(\yhom,\ad_{\mathbb C}\rhom))
\]
given by
\[
(x_1,\ldots,x_m) \mapsto \left(\sum\,\; x_i\,D\tilde h_i(\ahom),
\sum\,\;x_i\,\Hess\tilde h_i (\ahom)\right)
\]
is surjective. Here, $\Herm(V)$ stands for the Hermitian forms on
a complex vector space $V$. If we perturb $\tilde\Gamma$ by a small
amount, its projection $\Gamma = p(\tilde\Gamma)$ will be an embedded
link in $Y$.  It is clear that $\Gamma$ is mod--$2$ trivial; we claim
that in fact it is abundant.  In the discussion that follows, the
perturbing functions on $Y$ will be the push-down of the perturbing
functions $\tilde h_i$ on $\yhom$.  This makes sense because
the holonomy of $\ahom$ around a component $\tilde\gamma$ of
$\tilde\Gamma$ is the same as the holonomy of $A$ around
$p(\tilde\gamma)$.

According to Lemma \ref{L1:4orbit}, we have the following commutative
diagram

\[
\begin{CD}
\mathbb R^m @>>> \Hom(H^1(\yhom;\mathbb R);\mathbb R)\,\oplus\,
\Herm(H^1(\yhom;\ad_{\mathbb C}\rhom)) \\
@| @VV{(p^*)^*}V\\
\mathbb R^m @>>> \Hom(V^+(\rho);\mathbb R)\,\oplus\,\Sym(V^-(\rho))
\end{CD}
\]

\medskip\noindent
where the horizontal arrows are the holonomy maps as in~\eqref{E:holmap}.
The arrow along the top is surjective. According to Lemma \ref{L1:4orbit},
both the map $V^+(\rho)\to H^1(\yhom;\mathbb R)$ and the map $V^-(\rho)
\to H^1(\yhom;\ad_{\mathbb C}\rhom)$ are injective. By Lemma
\ref{L2:4orbit}, the map
\[
V^-(\rho)\to H^1(\yhom;\ad_{\mathbb C}\rhom)
\]
is obtained by complexification. Therefore, the right arrow in the diagram
is surjective, and hence so is the arrow on the bottom.
\end{proof}


\section{Proof of Theorem \ref{T:main}}
Let $0\ne w \in H^2 (Y;\Z_2)$ and consider a $U(2)$--bundle $P$ with
$c_1 (P) = w \pmod 2$. If $w$ is not in the image of $\iota:
H^2 (\pi_1 Y; \Z_2) \to H^2 (Y; \Z_2)$ then the theorem follows from
Corollary \ref{C:iota}. Otherwise, choose a 2--cocycle $c$ so that
$[c] = w$ and identify $\M(P)$ with $\P_c(Y;SU(2))$.

If $\M(P)$ is non--degenerate then Theorem \ref{T:main} follows because
no orbit in $\P_c (Y;\allowbreak SU(2))$ consists of one element, see
Remark \ref{R:one-orbits}, the contribution of the two--orbits equals
$(a_1\cup a_2\cup a_3)\,[Y]\pmod 2$ according to Proposition
\ref{P:cup}, and the orbits consisting of four and eight elements do
not contribute to $\lambda'''(Y,w)\pmod 2$ at all.

In general, $\M(P)$ needs to be perturbed to make it non--degenerate.
The two--orbits are already non--degenerate and hence if our
perturbation $h$ is sufficiently small they will remain such. The
perturbation $h$ will not create orbits with one element or new orbits
with two elements. Moreover, one can always achieve the  non--degeneracy
by using perturbations which are invariant with respect to the action of
$H^1(Y;\Z_2)$, see Corollary \ref{C:one} and Corollary \ref{C:two}.
Therefore, the above argument discarding the orbits with more than two
elements can be applied again to complete the proof of Theorem
\ref{T:main}.


\section{The Casson and Rohlin invariants for integral homology spheres}
In this section we explain how our Theorem \ref{T:main} implies Casson's
original result that $\lambda(\Sigma)= \rho(\Sigma)\pmod 2$ for integral
homology spheres $\Sigma$.

\subsection{Calculating the Casson invariant}
Every integral homology sphere $\Sigma$ can be obtained from $S^3$ by
surgery on an algebraically split link, that is, a link $k_1\cup\ldots
\cup k_n$ such that $\lk (k_i,k_j) = 0$ for $i\ne j$. Moreover, all
the surgery coefficients can be chosen to be $1$ or $-1$, so that
\begin{equation}\label{E:Sigma}
\Sigma = S^3 + \ep_1\cdot k_1 + \ldots + \ep_n\cdot k_n, \quad
\ep_i = \pm 1.
\end{equation}
The Casson invariants of  $\Sigma$ and  $\Sigma\pm k$ are related by
Casson's surgery formula
\[
\lambda(\Sigma\pm k) = \lambda(\Sigma)\pm \lambda' (\Sigma + 0\cdot k)
\]
where $\Sigma + 0\cdot k$ is the result of $0$--surgery of $\Sigma$
along $k$. In Casson's original approach, the term $\lambda'$ was
interpreted in terms of the Alexander polynomial of the knot $k$. For
our purposes, we interpret it gauge-theoretically.

Namely, let $P$ be a $U(2)$ bundle over $\Sigma + 0\cdot k$ such that
$w_2(\bar P)$ is dual to $[k]\in H_1 (\Sigma + 0\cdot k; \Z_2)$. Then
$\lambda'(\Sigma + 0\cdot k)$ is half a signed count of projectively
flat connections in $P$ with a fixed central part, modulo the gauge
group consisting of automorphisms of $P$ with determinant one (perhaps
after a perturbation). These projectively flat connections are counted
with signs determined by the Floer index. Therefore, the invariant
$\lambda' (\Sigma + 0\cdot k)$ equals half the Euler characteristic of
the Floer homology $I_*(\Sigma + 0\cdot k)$ so that the surgery formula
(\ref{E:Sigma}) is just an application of the Floer exact triangle

\medskip

\begin{picture}(250,84)
       \put(152,64)          {$I_* (\Sigma + 0\cdot k)$}
       \put(128,45)          {$\empty_{Y_*}$}
       \put(152,50)          {\vector(-1,-1){20}}
       \put(233,45)          {$\empty_{X_*}$}
       \put(229,30)          {\vector(-1,1){20}}
       \put(110,12)          {$I_* (\Sigma)$}
       \put(177,17)          {$\empty_{Z_*}$}
       \put(160,12)          {\vector(1,0){50}}
       \put(230,12)          {$I_* (\Sigma - k)$}
\end{picture}

\medskip\noindent
The surgery formula evaluates $\lambda(\Sigma \pm k)$ in terms of
$\lambda(\Sigma)$. Surgering out one knot at a time in the surgery
presentation (\ref{E:Sigma}), we end up with $S^3$ whose Casson
invariant is known to be trivial.

In order to calculate $\lambda(\Sigma)$ using this approach we need
to know the invariants $\lambda'(\Sigma + 0\cdot k)$ at each of the
steps. To this end, we use another surgery formula
\begin{equation}\label{E:last}
\lambda' (\Sigma + 0\cdot k \pm \ell) = \lambda' (\Sigma + 0\cdot k)
\pm \lambda'' (\Sigma + 0\cdot k + 0\cdot\ell)
\end{equation}
where $k\cup\ell$ is an algebraically split link in $\Sigma$ (it is
sufficient to work with algebraically split links because such
is the link in presentation (\ref{E:Sigma})). The term $\lambda''$
here is defined exactly as $\lambda'$ with only difference that now
$P$ is a $U(2)$ bundle such that $w_2 (\bar P)$ is dual to $[k] +
[\ell] \in H_1 (\Sigma + 0\cdot k + 0\cdot\ell; \Z_2)$. Again, the
above surgery formula follows from the Floer exact triangle, see
\cite{braam-donaldson:knots}.

This reduces calculation of $\lambda(\Sigma)$ to that of the
$\lambda''$--invariants. Applying the surgery formula yet another
time, we reduce the latter calculation to identifying $\lambda'''
(Y,w)$ for the homology torus $Y$ obtained by $0$--surgery on an
algebraically split link $k\cup\ell\cup m$ in $\Sigma$ with $w$
dual to $[k] + [\ell] + [m] \in H_1 (Y; \Z_2)$. Theorem \ref{T:main}
tells us that $\lambda'''(Y,w)$ equals $(a_1\cup a_2\cup a_3)\,[Y]
\pmod 2$ for any choice of basis $a_1$, $a_2$, $a_3\in H^1(Y;\Z_2)$.

\begin{remark}
A caveat in the above argument is that the simplification scheme it
is based upon fails for computing $\lambda'(S^3 + 0\cdot k)$. After
we simplified $\Sigma$ to $S^3$, a new scheme is needed to simplify
the knot, not the manifold itself. Such a simplification scheme,
based on skein moves, can be found in \cite{akbulut-mccarthy} or
\cite{saveliev:casson}. Again, it reduces calculation of $\lambda'
(S^3 + 0\cdot k)$ to that of $\lambda'''(Y,w)$.
\end{remark}

\subsection{Calculating the Rohlin invariant}
To conclude that $\lambda(\Sigma) = \rho(\Sigma)\pmod 2$ for all
integral homology spheres $\Sigma$, we will show that the Rohlin
invariant satisfies the same surgery formulas as the Casson invariant,
only reduced modulo 2.  Results similar to those in this section were found earlier by Turaev~\cite{turaev:linking}.

\begin{lemma}\label{L:last}
Let $\rho'(\Sigma + 0\cdot k)$ be the sum, over the two spin structures
on $\Sigma + 0\cdot k$, of their Rohlin invariants.  Then $\rho (\Sigma
+ k) = \rho (\Sigma) + \rho'(\Sigma + 0\cdot k)\pmod 2$.
\end{lemma}

\begin{proof}
The manifold $\Sigma + 0\cdot k$ can be obtained by $0$--surgery on
both $\Sigma + k$ and $\Sigma$.  Let $W_1$ and $W_2$ be the traces
of these surgeries, that is, smooth 4--manifolds obtained from $[0,1]
\times (\Sigma + k)$, respectively, $[0,1]\times\Sigma$, by attaching
a 2--handle along $\{\,1\,\} \times k$ with zero framing.  Then $W_1$
is a spin cobordism between $\Sigma + k$ and $\Sigma + 0\cdot k$ with
one spin structure, and $W_2$ is a spin cobordism between $\Sigma$
and $\Sigma + 0\cdot k$ with the other spin structure. Since the
intersection forms of both $W_1$ and $W_2$ are zero, we are finished.
\end{proof}

Note that changing the surgery coefficient from plus to minus does
not affect the Rohlin invariant, therefore, we may assume for the
sake of simplicity that all the surgery coefficients $\ep_i$ in
(\ref{E:Sigma}) are equal to one.

Let $k \cup \ell$ be an algebraically split link in $\Sigma$ and
define $\rho''(\Sigma + 0\cdot k + 0\cdot\ell)$ as the sum, over the
four spin structures on $\Sigma + 0\cdot k + 0\cdot\ell$, of their
Rohlin invariants. An argument similar to that of Lemma \ref{L:last}
proves the surgery formula
\[
\rho'(\Sigma + 0\cdot k + \ell) = \rho'(\Sigma + 0\cdot k) +
\rho'' (\Sigma + 0\cdot k + 0\cdot\ell),
\]
compare with (\ref{E:last}), and yet another application of the same
argument yields the formula
\[
\rho''(\Sigma + 0\cdot k + 0\cdot\ell) = \rho (\Sigma + k + \ell) +
\rho(\Sigma + k) + \rho(\Sigma + \ell) + \rho(\Sigma).
\]

This reduces the calculation of $\rho(\Sigma)$ to that of the
$\rho''$--invariants. Applying the surgery formula one more time, we
reduce the latter calculation to identifying $\rho'''(Y)$ for a
homology torus $Y = \Sigma + 0 \cdot k + 0\cdot\ell + 0\cdot m$. An
argument similar to that of Lemma \ref{L:last} yields
\begin{multline}\notag
\rho'''(Y) = \rho(\Sigma + k + \ell + m) + \rho(\Sigma + \ell + m) +
\rho(\Sigma + k + m) \\ + \rho(\Sigma + k + \ell) + \rho(\Sigma + k)+
\rho(\Sigma + \ell) + \rho(\Sigma + m) + \rho(\Sigma),
\end{multline}
which equals $(a_1 \cup a_2 \cup a_3)\,[Y]\pmod 2$ for any choice of
basis $a_1$, $a_2$, $a_3\in H^1 (Y; \Z_2)$, see 
\cite[Lemma 6.3]{kaplan:even}. This proves that $\lambda'''(Y,w) = 
\rho'''(Y)\pmod 2$ and therefore completes the proof of the formula 
$\lambda(\Sigma) = \rho(\Sigma)\pmod 2$.

\bigskip
\bigskip

\providecommand{\bysame}{\leavevmode\hbox to3em{\hrulefill}\thinspace}

\end{document}